\documentclass[reqno,12pt]{amsart}
\usepackage{amsmath,amssymb,amsthm}
\usepackage{a4wide}
\usepackage{hyperref}
\hypersetup{colorlinks=true,linkcolor=blue,citecolor=green,urlcolor=blue}
\numberwithin{equation}{section}
\allowdisplaybreaks[4]

\newtheorem{lemma}{Lemma}[section]
\newtheorem{theorem}{Theorem}[section]
\newtheorem{definition}{Definition}[section]
\newtheorem{proposition}{Proposition}[section]
\newtheorem{remark}{Remark}[section]
\newtheorem{corollary}{Corollary}[section]

\newcommand{\dis}{\displaystyle}

\newcommand{\R}{\mathbb{R}}

\renewcommand{\S}{\mathbb{S}}

\newcommand{\ga}{\gamma}
\newcommand{\om}{\omega}

\newcommand{\ta}{\theta}

\newcommand{\eqdef}{\overset{\mbox{\tiny{def}}}{=}}

\newcommand{\norm}[1]{\left\Vert #1 \right\Vert}
\newcommand{\moc}[1]{\left\lvert #1 \right\rvert}

\theoremstyle{plain}

\begin{document}
	\title[Boltzmann equations with dilatation and shear]{Boltzmann equations with dilatation and shear:\\ global well-posedness and stability}
\author[S.-Q. Liu]{Shuangqian Liu}
\address[SQL]{School of Mathematics and Statistics, and Key Lab NAA-MOE, Central China Normal University, Wuhan 430079, P.R.~China}
\email{sqliu@ccnu.edu.cn}

\author[H.-R. Peng]{Haoran Peng}
\address[HRP]{School of Mathematics and Statistics, and Hubei Key Lab--Math. Sci., Central China Normal University, Wuhan 430079, P.R.~China}
\email{mathematicphr@mails.ccnu.edu.cn}

\author[J.-T. Zou]{Jiatong Zou}
\address[JTZ]{School of Mathematics and Statistics, and Hubei Key Lab--Math. Sci., Central China Normal University, Wuhan 430079, P.R.~China}
\email{jiatongzou000818@mails.ccnu.edu.cn}

\date{\today}

\renewcommand{\subjclassname}{2020 Mathematics Subject Classification}	


	\begin{abstract}
We study global well-posedness and stability for Boltzmann equations
with linear velocity deformation, allowing for arbitrarily large deformation strength.
For diagonal dilatation with $-2<\gamma\leq0$ and cylindrical dilatation
with $-3/2<\gamma\leq0$, with or without shear, we construct global
homogeneous solutions for small data and prove their stability under
small spatially inhomogeneous perturbations. In scaled variables, the
perturbations remain uniformly controlled in deformation-adapted spaces
and decay in mixed spatial--velocity norms, algebraically in the diagonal
regime and exponentially in the cylindrical regime. For homogeneous
simple shear with $-2\leq\gamma\leq0$, we prove global existence and
uniqueness for nonnegative data with finite mass and Fisher information,
without smallness or a finite-energy assumption. These solutions conserve
mass and propagate additional Sobolev regularity. For spatially
inhomogeneous simple shear with $-2<\gamma\leq0$, we establish global
well-posedness for small data. The analysis combines commuting vector
fields, dispersive and bilinear collision estimates, and integrability
of the rescaled collision coefficients. For homogeneous shear, a
positivity-preserving splitting transfers Fisher-information control to
rough data, while a Lorentz-space estimate includes the endpoint
$\gamma=-2$ in the global continuation argument.
\end{abstract}

	\maketitle

\tableofcontents

\section{Introduction}\label{sec:introduction}
\subsection{The problem}
The Boltzmann equation describes the motion of a dilute gas through free
transport and binary collisions. In three space dimensions it reads
\begin{equation}\label{eq:intro-boltzmann}
 \partial_t f+v\cdot\nabla_x f=Q(f,f),
 \qquad (x,v)\in\mathbb R^3\times\mathbb R^3.
\end{equation}
Here, the unknown function $f(t,x,v)$ is the density distribution of the dilute gas at time $t\geq0$, position $x\in\R^3$ and velocity $v\in\R^3.$
Throughout this paper, the bilinear collision operator is defined by
\begin{align}\label{eq:collision-definition}
 Q(f,g)(v)=&\int_{\mathbb R^3}\int_{\mathbb S^2}\tilde{B}(v-v_*,\omega)
 \{f(v')g(v_*')-f(v)g(v_*)\}\,d\omega\,dv_*\notag\\
 \eqdef& Q^+(f,g)-Q^-(f,g),
\end{align}
where
\[
 v'=v-[(v-v_*)\cdot\omega]\omega,
 \qquad v_*'=v_*+[(v-v_*)\cdot\omega]\omega.
\]
Thus collisions conserve momentum and kinetic energy. We assume the
following form of Grad's angular cutoff condition:
	\begin{eqnarray}\label{eq:angular-cutoff}
		&\dis \tilde{B}(v-v_\ast,\om)=|v-v_\ast|^\ga B_0(\cos\ta),\ \cos\ta=\om\cdot\frac{v-v_\ast}{|v-v_\ast|},\ \om\in\S^2,\\
		&0\leq B_0(\cos\ta)\leq C|\cos\ta|.\notag
	\end{eqnarray}
Moreover, the collision frequency index $\ga\in(-3,1]$, characterizes the potential type, with $\ga\in[0,1]$ corresponding to hard potentials and $\ga\in(-3,0)$ corresponding to soft potentials. See
\cite{CercignaniIllnerPulvirenti1994} for the kinetic background.

An important class of nonequilibrium motions is given by homoenergetic
flows; see \cite{TruesdellMuncaster1980}. Let $A_0$ be a constant matrix such that
$\det(I+tA_0)>0$ for $t\geq0$, and set
\begin{equation*}
 L(t)=A_0(I+tA_0)^{-1}.
\end{equation*}
Then $L'+L^2=0$. The ansatz
$f(t,x,v)=g(t,v-L(t)x)$ reduces
\eqref{eq:intro-boltzmann} to
\begin{equation}\label{eq:intro-homo-equation}
 \partial_t g-L(t)v\cdot\nabla_v g=Q(g,g).
\end{equation}
Here the velocity variable is measured relative to the affine bulk
motion $L(t)x$. In particular, if the mass is finite, then
\[
 \int_{\mathbb R^3}g(t,v)\,dv
 =\frac{\int_{\mathbb R^3}g(0,v)\,dv}{\det(I+tA_0)}.
\]
The geometry of $L$ determines the competition between deformation
and collisions. The classification in \cite{JNVProfiles} includes
homogeneous dilatation, cylindrical dilatation with or without shear,
and simple shear, as well as other combinations of dilatation and shear.
We concentrate on the first three types of motion.

To study spatially dependent dynamics under the same velocity
deformation, we consider the model
\begin{equation}\label{eq:intro-inhomogeneous-model}
 \left\{
 \begin{aligned}
 &\partial_t f+v\cdot\nabla_x f-L(t)v\cdot\nabla_v f=Q(f,f),\\
 &f(0,x,v)=f_0(x,v).
 \end{aligned}\right.
\end{equation}
The deformation matrix $L(t)$ depends only on time. When $f$ is
independent of $x$, the spatial transport term vanishes and
\eqref{eq:intro-homo-equation} is recovered. Our spatially inhomogeneous
results below are formulated for \eqref{eq:intro-inhomogeneous-model}.

In this paper, we aim to construct small global homogeneous
solutions in the dilatational regimes and establish their stability under small spatially inhomogeneous
perturbations. For homogeneous simple shear, we investigate global existence
and uniqueness for nonnegative data with finite mass and Fisher information,
without a smallness condition. For the corresponding equation with
spatial transport, we study the global well-posedness for small
spatially inhomogeneous data near vacuum.

\subsection{Dilatational scaling and the perturbation equation}
We write the dilatational matrices in the common form
\begin{equation}\label{eq:dilatation-hypothesis}
 L(t)=\frac{P}{1+t}+R(t),\qquad
 R\in C([0,\infty);\mathbb R^{3\times3}),\qquad
 \|R(t)\|\leq C_R(1+t)^{-2},
\end{equation}
where either
\begin{equation}
 P=I\quad\hbox{or}\quad
 P=P_K:=\begin{pmatrix}1&0&K\\0&1&0\\0&0&0\end{pmatrix},
 \qquad K\in\mathbb R.\notag
\end{equation}
The first choice gives homogeneous diagonal dilatation; the second
gives cylindrical dilatation, with shear when $K\ne0$.
No smallness of $K$ or $C_R$ is imposed. The hypotheses
\eqref{eq:dilatation-hypothesis} are used directly in the analysis;
the remainder need not commute with $P$.

Put $a=\operatorname{tr}P/3$, so that $a=1$ or $a=2/3$, and define
\begin{equation}\label{eq:dilatation-scaling}
 \tau=\log(1+t),\qquad \xi=(1+t)^a v,\qquad
 m(t)=\exp\left(-\int_0^t\operatorname{tr}R(s)\,ds\right).
\end{equation}
With $f(t,x,v)=m(t)F(\tau,x,\xi)$ and
$g(t,v)=m(t)G(\tau,\xi)$, the scaled equations become
\begin{equation}\label{eq:scaled-system}
 \begin{aligned}
 \partial_\tau G-\nabla_\xi\cdot(A(\tau)\xi G)
   &=k(\tau)Q(G,G),\\
 \partial_\tau F+e^{(1-a)\tau}\xi\cdot\nabla_xF
       -\nabla_\xi\cdot(A(\tau)\xi F)
   &=k(\tau)Q(F,F),
 \end{aligned}
\end{equation}
where
\begin{equation}\label{eq:scaled-coefficients}
 \begin{aligned}
 A(\tau)&=P-aI+r(\tau),&
 r(\tau)&=e^\tau R(e^\tau-1),\\
 k(\tau)&=M(\tau)e^{-[a(3+\gamma)-1]\tau},&
 M(\tau)&=m(e^\tau-1).
 \end{aligned}
\end{equation}
In particular, $\|r(\tau)\|\leq C_Re^{-\tau}$ and $M$ is bounded
above and below by positive constants. Consequently,
\begin{equation}\label{eq:collision-time-integrability}
 k\in L^1(0,\infty)\quad\hbox{if}\quad
 \begin{cases}
 \gamma>-2,&P=I,\\
 \gamma>-3/2,&P=P_K.
 \end{cases}
\end{equation}
These are the two parameter ranges used for dilatation.

To study stability, write
\begin{equation*}
 F=G+H,\qquad F_0=G_0+H_0.
\end{equation*}
Then the perturbation satisfies
\begin{equation}\label{eq:intro-perturbation-equation}
 \begin{aligned}
 \partial_\tau H+e^{(1-a)\tau}\xi\cdot\nabla_xH
       -\nabla_\xi\cdot(A(\tau)\xi H)
 &=k(\tau)\{Q(H,H)+Q(H,G)+Q(G,H)\},\\
 H(0,x,\xi)&=H_0(x,\xi).
 \end{aligned}
\end{equation}
Since a nonzero homogeneous background does not belong to
$L^p(\mathbb R^3_x\times\mathbb R^3_\xi)$ for finite $p$, we estimate
$G$ in velocity spaces and $H$ in spatial--velocity spaces. The total
solution belongs to the resulting affine class.

\subsection{Main results}
We first specify the norms appearing in the statements. When
$-2<\gamma<0$, fix
\begin{equation}\label{eq:exponents}
 p_c=\frac{3}{3+\gamma},\qquad
 1<p_-<p_c<p_+<3,\qquad
 \frac1{p_-}+\frac1{p_+}=\frac2{p_c}.
\end{equation}
We use $E_T(h)$ for the inhomogeneous norm and $\mathcal H_T(g)$
for the homogeneous norm defined in Section~\ref{sec-nfs}.
Their initial-time versions are equivalent to
\begin{equation*}
 \begin{aligned}
 E_0(h)&\thicksim
 \begin{cases}
 \|h\|_{W^{6,1}_{x,\xi}},&\gamma=0,\\
 \|h\|_{W^{8,p_+}_{x,\xi}}+\|h\|_{W^{11,p_-}_{x,\xi}},&\gamma<0,
 \end{cases}\\
 \mathcal H_0(g)&\thicksim
 \begin{cases}
 \|g\|_{W^{6,1}_\xi},&\gamma=0,\\
 \|g\|_{W^{11,p_+}_\xi}+\|g\|_{W^{11,p_-}_\xi},&\gamma<0.
 \end{cases}
 \end{aligned}
\end{equation*}
Here $W^{m,p}$ is equipped with the sum of the $L^p$ norms of all
weak derivatives of order at most $m$. The time-dependent derivatives
in $E_T$ and $\mathcal H_T$ follow the deformation; the initial spaces
are ordinary, unweighted Sobolev spaces.

\begin{theorem}[Global dilatational solutions and their stability]
\label{thm:dilatation}
Assume \eqref{eq:angular-cutoff} and \eqref{eq:dilatation-hypothesis} hold. Let either $P=I$ and
$-2<\gamma\leq0$, or $P=P_K$ and $-3/2<\gamma\leq0$.
There are constants $\varepsilon_0>0$ and $C>0$, depending only
on the fixed kernel, exponents, and deformation, with the following
property. If
\[
 G_0\geq0,\qquad F_0=G_0+H_0\geq0,\qquad
 \mathcal H_0(G_0)+E_0(H_0)\leq\varepsilon_0,
\]
then \eqref{eq:scaled-system} has a unique global mild homogeneous
solution $G$ and a unique global mild solution $F=G+H$ in the class
$\mathcal H_\infty(G)+E_\infty(H)<\infty$. They are nonnegative and
satisfy
\begin{equation*}
 \mathcal H_\infty(G)\leq C\mathcal H_0(G_0),\qquad
 E_\infty(H)\leq CE_0(H_0).
\end{equation*}
Moreover, for every $\tau\geq0$,
\begin{equation*}
 \sum_{|\alpha|+|\beta|\leq m}
 \|\partial_x^\alpha Y^\beta H(\tau)\|_{L_x^\infty L_\xi^p}
 \leq C [d(\tau)]^{-1/p}E_0(H_0),
 \qquad
 d(\tau)=
 \begin{cases}
 (1+\tau)^3,&P=I,\\
 e^\tau(1+\tau)^2,&P=P_K.
 \end{cases}
\end{equation*}
Here $(p,m)=(1,3)$ when $\gamma=0$, and
$(p,m)\in\{(p_-,8),(p_+,5)\}$ when $-2<\gamma<0$.
The commuting fields $Y$ are defined in Section~\ref{sec-nfs}.
The map $G_0\mapsto G$ is locally Lipschitz in the stated norms;
for a fixed background $G$, the same holds for $H_0\mapsto H$.
If also $G_0\in L^1_\xi$, then
$\int G(\tau,\xi)\,d\xi=\int G_0(\xi)\,d\xi$.
\end{theorem}

We now turn to simple shear with unit coefficient. The homogeneous problem is
\begin{equation}\label{sim-shear-hom}
 \partial_tg-v_2\partial_{v_1}g=Q(g,g),\qquad g(0,v)=g_0(v),
\end{equation}
and the spatially inhomogeneous problem is
\begin{equation}\label{eq:shear-cauchy}
 \left\{
 \begin{aligned}
 &\partial_t f+v\cdot\nabla_xf-v_2\partial_{v_1}f=Q(f,f),\\
 &f(0,x,v)=f_0(x,v).
 \end{aligned}\right.
\end{equation}
Now we state the homogeneous result.

\begin{theorem}[Global well-posedness under simple shear]
\label{thm:shear}
Under assumption \eqref{eq:angular-cutoff}, with notation in Definition \ref{defi-norm} and \ref{defi-shear-hom}, the following assertions hold.
\begin{enumerate}
\renewcommand{\labelenumi}{\textup{(\roman{enumi})}}
\renewcommand{\theenumi}{\roman{enumi}}
\item\label{part:shear-homogeneous}
Let $-2\leq\gamma\leq 0$ and let $g_0\geq0$ satisfy
\[
 M_0:=\int_{\mathbb R^3}g_0(v)\,dv<\infty,
 \qquad \mathcal I(g_0)<\infty.
\]
Then the Cauchy problem \eqref{sim-shear-hom} admits a unique nonnegative global mild solution $g\in C([0,\infty),Z)$,
satisfying, for all $t\geq0$,
\begin{equation*}
 \int_{\mathbb R^3}g(t,v)\,dv=M_0,\qquad
 \mathcal I(g(t))\leq e^t\mathcal I(g_0).
\end{equation*}
If, in addition,
$g_0\in Z^m$ for an integer $m\geq0$, then $g\in C([0,\infty);Z^m)$.

\item\label{part:shear-inhomogeneous}
Let $-2<\gamma\leq0$.
There are constants $\varepsilon_1>0$ and $C>0$ such that, for every
$f_0\geq0$ satisfying $E_0(f_0)\leq\varepsilon_1$,
the Cauchy problem \eqref{eq:shear-cauchy} admits a unique nonnegative global mild solution
with
\begin{equation*}
 E_\infty(f)\leq CE_0(f_0).
\end{equation*}
Here $E_T$ uses the simple-shear vector fields of
Section~\ref{sec-nfs}. The data-to-solution map is locally Lipschitz
from the initial space to this global solution space.
\end{enumerate}
\end{theorem}

\begin{remark}\label{rem:scope}
In Theorem~\ref{thm:dilatation}, decay is measured in scaled velocity
and time $\tau=\log(1+t)$; exponential decay in $\tau$ is therefore
algebraic in the original time.
Theorem~\ref{thm:shear}\textup{(\ref{part:shear-homogeneous})}
requires neither small mass nor finite energy. Its Fisher bound
ensures finite-time continuation, without asserting convergence to
equilibrium. Theorem~\ref{thm:shear}\textup{(\ref{part:shear-inhomogeneous})}
concerns small spatially inhomogeneous data near vacuum on $\mathbb R^3_x$.
\end{remark}

\subsection{Related works}
We briefly review homoenergetic and shear flows, near-vacuum dynamics,
and Fisher-information estimates relevant to the present work.

\subsubsection{Homoenergetic solutions and shear flows}
Early studies of homoenergetic flows include the work of Galkin
\cite{G1} and Truesdell \cite{Truesdell1956} on kinetic moment equations;
see also Truesdell and Muncaster \cite{TruesdellMuncaster1980}.
Cercignani \cite{Cercignani1989} developed an existence theory for
homoenergetic affine flows, while Bobylev, Caraffini and Spiga
\cite{BobylevCaraffiniSpiga1996} studied the associated group-invariant
reductions. For homogeneous Maxwell models without deformation,
self-similar solutions and their asymptotics, including infinite-energy
solutions, were studied in \cite{BobylevCercignani2002,BobylevCercignani2003};
a Fourier-based measure theory was developed in \cite{CannoneKarch2010}.

James, Nota and Vel\'azquez distinguished balanced, collision-dominated
and hyperbolic-dominated regimes, constructing self-similar profiles
\cite{JNVProfiles} and deriving formal asymptotic descriptions
\cite{JNVCollision,JNVHyperbolic}. For Maxwell molecules with small
linear deformation, profile stability and self-similar asymptotics
were established in \cite{BobylevNotaVelazquez2020,Bobylev2023}.
Kepka \cite{Kepka2023} treated noncutoff Maxwell molecules with small
linear deformation, including
smooth profiles and stability with moments of order greater than two.
Anisotropic rescaling for objective Boltzmann and Fokker--Planck
solutions was studied by Matthies and Theil \cite{MatthiesTheil2019}.

For cutoff uniform shear, Duan and Liu \cite{DuanLiuUSF} constructed
self-similar Maxwell-molecule profiles and proved stability under
spatially inhomogeneous perturbations for small shear rates.
Their hard-potential analysis \cite{DuanLiuUSFHard} gives asymptotic
expansions and algebraic temperature growth for small deformation and
suitably prepared data. Kepka \cite{Kepka2024} also obtained Maxwellian
asymptotics for cutoff and noncutoff hard potentials in the
collision-dominated regime, using a different approach. For cutoff Maxwell molecules and hard
potentials with an additional thermostat, small deformation admits
non-Maxwellian steady states stable under periodic perturbations \cite{DuanLiuThermostat2022}.
Recently, Liu, Wang, Yang and Zhang \cite{LiuWangYangZhang2025} proved perturbative stability for homogeneous hard potentials under small constant deformations with negative trace.

For boundary-driven flows, Esposito, Lebowitz and Marra
\cite{EspositoLebowitzMarra1994,EspositoLebowitzMarra1995} established
stationary hydrodynamic limits in a slab under small forcing and
boundary perturbations. At a fixed kinetic scale, Duan, Liu and Yang
\cite{DuanLiuYangCouette2025} constructed stable plane Couette flows
for cutoff Maxwell molecules, small shear rates and diffuse reflection.
For spatially periodic uniform shear with Maxwell molecules,
Duan, Liu and Shen
\cite{DuanLiuShen2025} proved global existence and analyzed energy growth.
Recent complementary results include formal strong-shear asymptotics
\cite{NotaVelazquezStrongShear} and stationary states for the linear
Rayleigh--Boltzmann equation coupled to a prescribed bath
\cite{MieleNotaVelazquez2025}. The present paper concerns
whole-space stability near small homogeneous states under dilatation,
small-data existence under shear, and homogeneous shear continuation
from finite mass and Fisher information.

\subsubsection{Near-vacuum dynamics and related Cauchy theories}
Global renormalized existence for large data under integrability and entropy bounds \cite{DiPernaLions1989} does not supply general uniqueness. Near Maxwellians, collision dissipation supports classical solutions, including periodic cutoff soft potentials $-3<\gamma<0$ \cite{GuoCutoff2003}, and macro--micro energy methods \cite{LiuYangYu2004}. Anisotropic coercivity underlies noncutoff stability \cite{GressmanStrain2011}. Critical Besov and low-spatial-regularity theories treat cutoff \cite{DuanLiuXu2016,DuanSakamoto2018} and noncutoff kernels \cite{MorimotoSakamoto2016,DuanLiuSakamotoStrain2021}. These near-Maxwellian methods rely on the dissipative structure supplied by a nonzero equilibrium and cannot be applied directly to the vacuum problem, since the linearization at vacuum has no nontrivial collision dissipation. Homogeneous-background stability also predates the present work: small spatial variation is considered in \cite{ArkerydEspositoPulvirenti1987}, and small perturbations of radially symmetric homogeneous data with specular reflection in \cite{GuoLiu2017}.

Near vacuum, transport dispersion replaces linearized collision dissipation. The characteristic approach of Illner and Shinbrot \cite{IllnerShinbrot1984} extends to small-data Vlasov--Poisson--Boltzmann (VPB) problems \cite{GuoVacuum2001,DuanYangZhu2006} and external forces satisfying structural characteristic assumptions \cite{DuanYangZhu2005}. Collision potentials yield asymptotic comparison with linear Vlasov dynamics \cite{ChaeHaHwang2006}. For VPB, uniform $L^1$ stability in \cite{DuanZhangZhu2006} requires $N\geq4$; its $N=3$ result is local in time. Uniform $L^p$ stability for the unforced Boltzmann equation is established in \cite{HaYamazakiYun2008}.

Dispersive collision estimates give small-data mild solutions beyond weighted supremum spaces \cite{Arsenio2011}. For the two-dimensional constant-kernel equation, bilinear spacetime estimates yield global well-posedness and scattering for regular, localized data with small critical $L^2$ norm \cite{ChenDenlingerPavlovic2021}. Large scattering solutions near joint space--velocity Gaussians are constructed under additional physical constraints in \cite{ChenDenlingerPavlovic2024}, with uniqueness for Schwartz data.

In three dimensions, He and Jiang \cite{HeJiang2023} proved scattering for $-1\leq\gamma<0$ with small weighted $L^3$ data and additional integrability, correcting a loss-term estimate in \cite{HeJiang2017}. Sharp gain estimates and further scattering results appear in \cite{HeJiangKuoLiang2024}. Sobolev well/ill-posedness thresholds for constant kernels in three dimensions and soft kernels in two and three dimensions \cite{ChenHolmer2024,ChenShenZhangSoft} show that scaling alone does not determine well-posedness. Global scattering for $-\tfrac12\leq\gamma\leq0$ is proved under critical-norm smallness and additional Sobolev regularity \cite{ChenShenZhangScattering}. Recent advances include unconditional uniqueness through decoupling \cite{ChenShenZhangUniqueness2026} and anisotropic Besov estimates using a div--curl argument \cite{HuLiuPengZhou2026}.

For long-range interactions, weighted estimates and commuting fields underpin vacuum stability with moderately soft potentials for the Landau equation \cite{Luk2019} and the noncutoff Boltzmann equation \cite{Chaturvedi2021}. Small-data Vlasov--Yukawa--Boltzmann theories cover cutoff dynamics, stability and scattering \cite{ChoiHa2015}, and noncutoff moderately soft potentials \cite{HuXiongZhu2025}. These approaches motivate our estimates, while linear velocity deformation requires a separate analysis of anisotropic characteristics and background--perturbation interactions.

\subsubsection{Fisher information and homogeneous continuation}
Fisher-information estimates for Maxwell collisions were established
by Villani \cite{VillaniFisher1998}; Matthes and Toscani
\cite{MatthesToscani2012} gave a Fourier-based proof for the gain
operator. Guillen and Silvestre \cite{GuillenSilvestre2025} proved
Fisher-information monotonicity for a class of homogeneous Landau equations,
including the Coulomb case, and used it to obtain global regularity.
For Boltzmann collisions, the monotonicity theorem of Imbert, Silvestre
and Villani \cite{ISV} supplies the estimate used here. In the current paper, for the cutoff shear problem, combining monotonicity
with the drift and mass conservation yields $L^{3,1}_v$ control on every
finite time interval. This provides global continuation, including
$\gamma=-2$, assuming only finite mass and Fisher information.

The present paper addresses several issues that remain beyond the scope of those results: the stability of small homogeneous dilatational states under spatially inhomogeneous perturbations, the global near-vacuum dynamics under simple shear, and the continuation of homogeneous shear solutions assuming only finite mass and Fisher information, including the borderline case $(\gamma=-2)$. These problems provide the main motivation for the present work.

Our aim is to develop estimates capable of handling anisotropic deformations without imposing a smallness assumption on the deformation strength, as well as a continuation mechanism for homogeneous shear flows that requires neither small initial data nor finite energy. The distinct roles played by spatial dispersion and Fisher information naturally lead to the different analytical approaches developed below.


\subsection{Difficulties and main ideas}
The proofs rely on three ingredients, each addressing a different
obstruction created by the deformation.

\medskip
\noindent$\bullet$ \textit{Commuting fields and unweighted collision estimates.}
A fixed polynomial velocity weight grows along shear characteristics.
For example, if
$T_K=\partial_t+v\cdot\nabla_x-Kv_2\partial_{v_1}$ and
$w_\ell(v)=\langle v\rangle^\ell$, then
\[
 T_K(w_\ell f)=w_\ell Q(f,f)
       -K\ell\frac{v_1v_2}{\langle v\rangle^2}w_\ell f.
\]
The additional term has no favorable sign. We incorporate the full
linear deformation into the derivatives, following the commuting-field
approach to dispersion in \cite{Wong2018}. For the scaled
transport, solve
\begin{equation}\label{B-C-defi}
     \Phi'=-A\Phi,\qquad \Psi'=e^{(1-a)\tau}\Phi,\qquad \Phi(0)=\Psi(0)=I,
\end{equation}
and set $Y=\Phi^T\nabla_\xi+\Psi^T\nabla_x$. These fields commute with
the inhomogeneous linear operator, while $\Phi^T\nabla_\xi$ commutes
with its homogeneous counterpart. Translation covariance of $Q$
gives the corresponding Leibniz rule. The nonlinear estimates
therefore close in unweighted norms built from the same derivatives.
For soft potentials, exponents on both sides of $p_c$ control the
singular and nonsingular parts of $|\cdot|^\gamma*g$. The derivative
levels $8$ and $11$ allow the spatial Sobolev embedding needed when
these two velocity exponents interact. For simple shear,
$\det \Psi(t)=(1+t)^3$, and the resulting factor $(1+t)^{-3/p_+}$
is integrable because $p_+<3$. This explains the restriction
$\gamma>-2$ in the small-data shear argument.

\medskip
\noindent$\bullet$ \textit{Anisotropic volume growth and stability of an affine class.}
In the dilatational problem, homogeneous backgrounds do not disperse
in space. Thus the mixed terms $Q(H,G)$ and $Q(G,H)$ cannot be
estimated by the same decay argument used for $Q(H,H)$. We instead
use the integrable collision coefficient $k$ to control all three
terms in \eqref{eq:intro-perturbation-equation}. Quantitative
convergence is then obtained from the spatial part of the commuting
fields. The delicate geometric step is the cylindrical case:
$P_K$ is a projection, but it need not commute with the decaying
remainder $r$. An equation for the normalized matrix
$e^{-\tau/3}\Psi\Phi^{-1}$ separates the $P_K$ and $I-P_K$ components
without discarding the remainder. It yields
$\det \Psi \thicksim e^\tau(1+\tau)^2$ for sufficiently large $\tau$. This links the anisotropic
volume growth of the deformed flow to mixed-norm decay, while
preserving a separate velocity-space norm for the background.

\medskip
\noindent$\bullet$ \textit{Fisher control and global continuation without smallness.}
Homogeneous simple shear has neither spatial dispersion nor an
integrable prefactor in front of $Q$. We combine collision-induced
Fisher-information dissipation from \cite{ISV} with the growth
bound for the shear group. A positivity-preserving approximation
justifies this combination for data with finite mass and Fisher
information. The Sobolev inequality gives
\[
 \|g(t)\|_{L_v^3}=\|\sqrt{g(t)}\|_{L_v^6}^2
 \leq C\mathcal I(g(t))\leq Ce^t\mathcal I(g_0).
\]
Since $|v|^{-2}\in L^{3/2,\infty}$, this estimate and mass conservation
bound the collision frequency on every finite time interval for
$-2\leq\gamma\leq0$. Together with the gain estimate, this closes the collision operator in
$Z=L^1\cap L^{3,1}$, yielding global continuation without smallness or
an energy assumption. Fisher monotonicity and the Sobolev--Lorentz
inequality are inputs; the contribution here is their combination with
the shear evolution and the approximation argument at $\gamma=-2$.

\subsection{Organization of the paper}
Section~\ref{sec-nfs} introduces the commuting fields and the norms
used throughout the paper. Section~\ref{dis-esti} proves the common
dispersive and collision estimates. Section~\ref{sec-dial-cylind}
establishes the matrix bounds and proves Theorem~\ref{thm:dilatation}.
Section~\ref{sec-5} proves Theorem~\ref{thm:shear}, treating the
homogeneous case first and the inhomogeneous case second. The appendix \ref{sec:appendix}
is devoted to the auxiliary Sobolev and velocity-space estimates.

\section{Commuting vector fields and function spaces}\label{sec-nfs}
This section fixes the notation used in the proofs. We introduce derivatives adapted to the linear deformation and use a single pair of energy norms for the homogeneous solution and its spatial perturbation. The same notation applies to all three regimes.

Write
\[
 \mathcal Lh=\partial_\tau h+(\mathcal T(\tau)\xi)\cdot\nabla_xh
       -\nabla_\xi\cdot(A(\tau)\xi h),
 \qquad
 \mathcal L_Hg=\partial_\tau g-\nabla_\xi\cdot(A(\tau)\xi g),
\]
where $\mathcal T(\tau)$ will be given as \eqref{com-eqn}.
Let the matrix-valued functions $\Phi(\tau)$ and $\Psi(\tau)$ solve
\begin{equation}\label{eq:commuting-matrices}
 \Phi'=-A\Phi,\qquad \Psi'=\mathcal T\Phi,\qquad \Phi(0)=\Psi(0)=I,
\end{equation}
and set
\begin{equation}\label{eq:commuting-field-definition}
 X=\nabla_x,\qquad Y=\Phi^T\nabla_\xi+\Psi^T\nabla_x,
 \qquad Y_H=\Phi^T\nabla_\xi.
\end{equation}
Direct differentiation gives
$[X_j,\mathcal L]=[Y_j,\mathcal L]=0, [Y_{H,j},\mathcal L_H]=0.$
All spatial and velocity derivatives in \eqref{eq:commuting-field-definition} commute with one another. Moreover, $Yg=Y_Hg$ for an $x$-independent function $g$.

For multi-indices $\alpha,\beta\in\mathbb N_0^3$, write
$\mathcal D^{\alpha,\beta}=X^\alpha Y^\beta$. For $m\in\mathbb N_0$ and $1\leq p<\infty$, define
\[
 \norm{h}_{m,p}(\tau)=\sum_{|\alpha|+|\beta|\leq m}
       \|\mathcal D^{\alpha,\beta}(\tau)h\|_{L^p_{x,\xi}},
 \qquad
 \norm{g}^H_{m,p}(\tau)=\sum_{|\beta|\leq m}
       \|Y_H^\beta(\tau)g\|_{L^p_\xi}.
\]
We also write $\mathcal D^I$ for $\mathcal D^{\alpha,\beta}$, with $I=(\alpha,\beta)\in\mathbb N_0^6$.

\begin{definition}\label{defi-norm}
For $-2<\gamma\leq0$, use the exponents fixed in \eqref{eq:exponents}. Define the instantaneous energies by
\[
 E_\tau(h)=
 \begin{cases}\norm{h}_{6,1}(\tau),&\gamma=0,\\
 \norm{h}_{8,p_+}(\tau)+\norm{h}_{11,p_-}(\tau),&-2<\gamma<0,
 \end{cases}
\]
\[
 \mathcal H_\tau(g)=
 \begin{cases}\norm{g}^H_{6,1}(\tau),&\gamma=0,\\
 \norm{g}^H_{11,p_+}(\tau)+\norm{g}^H_{11,p_-}(\tau),&-2<\gamma<0.
 \end{cases}
\]
For $0<T\leq\infty$, set
\[
 E_T(h)=\sup_{0\leq\tau\leq T}E_\tau(h(\tau)),
 \qquad \mathcal H_T(g)=\sup_{0\leq\tau\leq T}\mathcal H_\tau(g(\tau)).
\]
When $T=\infty$, the supremum is over finite $\tau\geq0$. The notation $E_0(h_0)$ and $\mathcal H_0(g_0)$ refers to the corresponding initial-data energies. We denote the resulting solution spaces by $X_T$ and $X_{H,T}$, and write $X=X_\infty$, $X_H=X_{H,\infty}$. Initial-data spaces are the completions of smooth compactly supported functions in these norms; solutions are continuous after pullback by the linear flow.
\end{definition}

The admissible exponents can be chosen explicitly. Choose
\[
 0<\delta<\min\left\{\dfrac{2}{3}+\dfrac{\gamma}{3},\,-\dfrac{\gamma}{3}\right\},
 \qquad \frac1{p_-}=1+\dfrac{\gamma}{3}+\delta,
 \qquad \frac1{p_+}=1+\dfrac{\gamma}{3}-\delta.
\]
The interval for $\delta$ is nonempty precisely when
$-2<\gamma<0$, and these choices satisfy \eqref{eq:exponents}.
The two sides of the threshold have different uses: $p_+>p_c$
controls the singularity of the collision kernel at zero relative
velocity, while $p_-<p_c$ controls its tail without a polynomial
velocity weight. The additional restriction $p_+<3$ provides the necessary integrability to close the inhomogeneous simple-shear estimate. At $\gamma=-2$, that threshold reaches $3$,
and this integrability argument no longer applies. The cylindrical
range has the further restriction $\gamma>-3/2$ because the
homogeneous rescaled collision coefficient must also be integrable.
For $\gamma=0$, the two-exponent construction degenerates; the
$L^1$ collision estimate gives the separate endpoint theory used
in the definition.

We next introduce the velocity spaces and notation used to study
\eqref{sim-shear-hom} for $-2\leq\gamma\leq0$.
\begin{definition}\label{defi-shear-hom}
For an integer $m\geq0$, define
\[
 Z^m=W^{m,1}(\mathbb R^3)\cap W^{m,(3,1)}(\mathbb R^3),
 \qquad \|h\|_{Z^m}=\|h\|_{W^{m,1}}+\|h\|_{W^{m,(3,1)}},
 \qquad Z=Z^0.
\]
For $0<p<\infty$ and $0<q\leq\infty$, the Lorentz space $L^{p,q}$ is defined by the finiteness of the quasi-norm
\[\norm{f}_{L^{p,q}}=\begin{cases}
    \left(\displaystyle\int_0^\infty
    (t^{\frac 1p} f^*(t))^q \frac{dt}{t} \right)^{\frac 1q} & q<\infty, \\[4pt]
    \sup_{t>0} t^{\frac 1p} f^*(t) & q=\infty,
\end{cases}\]
where the decreasing rearrangement $f^*(t):=\inf\{\lambda\geq 0: |\{v:|f(v)|>\lambda\}|\leq t\}$. And The Lorentz-Sobolev space $W^{m,(3,1)}(\mathbb R^3)$ consists of functions whose distributional derivatives of every order up to $m$ belong to $L^{3,1}$.
For nonnegative $h$, its Fisher information is
\[
 \mathcal I(h)=4\|\nabla\sqrt h\|_{L^2}^2,
\]
with value $+\infty$ if $\sqrt h\notin H^1$ when $h\in L^1$. For smooth positive $h$, this equals $\int |\nabla h|^2/h$.
The shear group is
\[
 (S(t)h)(v)=h(v_1+tv_2,v_2,v_3),\qquad t\in\mathbb R.
\]
It is a positive isometry on $L^1,L^{3,1}$ and $Z$. Density of $C_c^\infty$ in these spaces implies strongly continuity.
\end{definition}
Throughout the paper, $C$ denotes a positive constant that may change from line to line. The notation $a\lesssim b$ and $a\gtrsim b$ respectively means $a\leq Cb$ and $a\geq Cb$, and $a\thicksim b$ means $a\lesssim b$ and $b\lesssim a$; constants may depend on the fixed collision kernel, exponents, and deformation, but not on the time interval.

\section{Linear dispersion and collision estimates}\label{dis-esti}
In this section, we establish the estimates common to the proofs of the main theorems. The linear geometry enters only through the growth of $\det \Psi$. Translation invariance of the collision operator then allows the same commuting derivatives to control homogeneous, inhomogeneous, and mixed interactions.


The first result is concerned with dispersive properties from the commuting vector fields.
\begin{lemma}\label{thm:dis-esti}
Let $d:[0,\infty)\to[1,\infty)$ be locally bounded. Suppose that, for some $\tau_0\geq0$ and $c>0$,
\[
 |\det \Psi(\tau)|\geq c\,d(\tau),\qquad \tau\geq\tau_0.
\]
Then, for $1\leq p<\infty$,
\begin{equation}\label{Es-Disper}
 \|h\|_{L_x^\infty L_\xi^p}
 \lesssim [d(\tau)]^{-1/p}\norm{h}_{3,p}(\tau).
\end{equation}
\end{lemma}
\begin{proof}
 For $\tau \in [0,\tau_0]$, by Lemma \ref{lem:banach-valued-sobolev},
      \begin{align}\label{Dis-1}
        \norm{f}_{L_{x}^{\infty}L_\xi^{p}}&\lesssim\sum_{|I|\leq 3} \norm{X^{I}f}_{L_{x,\xi}^p}\lesssim \sum_{|I|\leq 3} \norm{\mathcal{D}^{I}f}_{L_{x,\xi}^p}.
      \end{align}
For $\tau\in [\tau_0,\infty),$ using velocity translation and the invertibility of $\Psi$ gives
      \begin{align}\label{Dis-2}
     \notag   \norm{f}_{L_{x}^{\infty}L_\xi^{p}}&= \norm{f(t,\Psi x,\Phi x+\xi)}_{L_{x}^{\infty}L_\xi^{p}}\lesssim \sum_{|I|\leq 3}\norm{X^{I}[f(t,\Psi x,\Phi x+\xi)]}_{L_{x,\xi}^{p}}\\
                                    &\thicksim \sum_{|I|\leq 3} \norm{Y^{I}f(t,\Psi x,\Phi x+\xi)}_{L_{x,\xi}^{p}}
                                          \lesssim \sum_{|I|\leq 3} |\text{det}\Psi(\tau)|^{-\frac{1}{p}}\norm{\mathcal{D}^{I}f(t,x,\xi)}_{L_{x,\xi}^{p}}\\
        \notag                                   &\lesssim [d(\tau)]^{-\frac{1}{p}}\sum_{|I|\leq 3}\norm{\mathcal{D}^I f}_{L_{x,\xi}^p}.
      \end{align}
      Therefore, by \eqref{Dis-1} and \eqref{Dis-2}, we deduce that \eqref{Es-Disper}. This ends the proof of Lemma \ref{thm:dis-esti}.

\end{proof}

We shall use
\begin{equation*}
 d(\tau)=
 \begin{cases}
 (1+\tau)^3,&\text{diagonal dilatation and simple shear},\\
 e^\tau(1+\tau)^2,&\text{cylindrical dilatation, with or without shear}.
 \end{cases}
\end{equation*}
The dilatational determinant estimates are proved in Section~\ref{sec-dial-cylind}. For simple shear,
\begin{equation*}
 \Phi=I-\tau A,\qquad \Psi=(1+\tau)I-\tfrac12\tau^2A,
\end{equation*}
where $A=(A_{ij})_{3\times 3}$ satisfies $A_{ij} = 0$ but $A_{12}=1$.
so $\det \Phi=1$ and $\det \Psi=(1+\tau)^3$.

The following result is devoted to the bilinear estimates in terms of $ E_\tau$ and $ H_\tau$ norms given by Definition \ref{defi-norm}.

\begin{lemma}\label{thm:bilinear-estimation}
Assume the collision-kernel hypotheses stated in the introduction and the hypotheses of Lemma~\ref{thm:dis-esti}. Set $p_*=p_+$ for $-2<\gamma<0$ and $p_*=1$ for $\gamma=0$. Then, at each fixed time,
\begin{align}
    E_\tau(Q(f,g))
    &\lesssim [d(\tau)]^{-1/p_*}E_\tau(f)E_\tau(g),
    \label{Estimate-Q-X}\\
    \mathcal H_\tau(Q(F,G))
    &\lesssim \mathcal H_\tau(F)\mathcal H_\tau(G),
    \label{Estimate-Q-X-Hom}\\
    E_\tau(Q(f,G))
    &\lesssim E_\tau(f)\mathcal H_\tau(G),
    \label{Estimate-Q-X-interaction-1}\\
    E_\tau(Q(F,g))
    &\lesssim \mathcal H_\tau(F)E_\tau(g),
    \label{Estimate-Q-X-interaction-2}
\end{align}
where $F$ and $G$ are independent of $x$.
\end{lemma}

\begin{proof}
It suffices to consider the case $\gamma\in(-2,0)$, since the case $\gamma=0$ is similar and simpler.

Recalling \eqref{eq:collision-definition}, we split the collision operator $Q(\cdot,\cdot)$ into its gain and loss parts,
\[
    Q(\cdot,\cdot)
    =
    Q^{+}(\cdot,\cdot)-Q^{-}(\cdot,\cdot),
\]
and estimate them separately in the three interaction regimes considered below. For notational convenience, we set
\[
    p\in\{p_{+},p_{-}\},
    \qquad
    N_{p_{+}}=8,
    \qquad
    N_{p_{-}}=11.
\]

\medskip
\noindent\textbf{Step 1: Inhomogeneous interactions.}

We first estimate the gain term. Whenever
$|I_1|+|I_2|\leq N_p,$
at least one of $|I_1|$ and $|I_2|$ is at most five. By Lemma~\ref{lem:collision-velocity} and Lemma~\ref{thm:dis-esti}, we obtain
\begin{align*}
    \sum_{|I|\leq N_p}
    \norm{\mathcal{D}^{I}Q^{+}(f,g)}_{L_{x,\xi}^{p}}
    &\lesssim
    \sum_{|I_1|+|I_2|\leq N_p}
    \norm{Q^{+}(\mathcal{D}^{I_1}f,\mathcal{D}^{I_2}g)}_{L_{x,\xi}^{p}}
    \\
    &\lesssim
    \sum_{\substack{|I_1|+|I_2|\leq N_{p}\\ |I_1|\leq 5}}
    \norm{\mathcal{D}^{I_1} f}_{L_{x}^{\infty}L_{\xi}^{p_c}}
    \norm{\mathcal{D}^{I_2}g}_{L_{x,\xi}^{p}}
    \\
    &\quad
    +
    \sum_{\substack{|I_1|+|I_2|\leq N_{p}\\ |I_2|\leq 5}}
    \norm{\mathcal{D}^{I_1} f}_{L_{x,\xi}^{p}}
    \norm{\mathcal{D}^{I_2}g}_{L_{x}^{\infty}L_{\xi}^{p_c}}
    \\
    &\lesssim
    [d(\tau)]^{-\frac{1}{p_c}}
    \left\{
    \sum_{\substack{|I_1|+|I_2|\leq N_{p}\\ |I_1|\leq 5,\ |J|\leq 3}}
    \norm{\mathcal{D}^{I_1+J} f}_{L_{x,\xi}^{p_c}}
    \norm{\mathcal{D}^{I_2}g}_{L_{x,\xi}^{p}}
    \right.
    \\
    &\qquad\qquad\qquad\left.
    +
    \sum_{\substack{|I_1|+|I_2|\leq N_{p}\\ |I_2|\leq 5,\ |J|\leq 3}}
    \norm{\mathcal{D}^{I_1} f}_{L_{x,\xi}^{p}}
    \norm{\mathcal{D}^{I_2+J}g}_{L_{x,\xi}^{p_c}}
    \right\}
    \\
    &\lesssim
    [d(\tau)]^{-\frac{1}{p_c}}
    E_\tau(f)E_{\tau}(g),
\end{align*}
where we have used the interpolation inequality
\[
    \norm{\mathcal{D}^{I} h}_{L_{x,\xi}^{p_c}}
    \lesssim
    \norm{\mathcal{D}^{I} h}_{L_{x,\xi}^{p_{-}}}
    +
    \norm{\mathcal{D}^{I} h}_{L_{x,\xi}^{p_+}}.
\]
Summing over $p=p_{+}$ and $p=p_{-}$ gives
\[
    E_{\tau}(Q^{+}(f,g))
    \lesssim
    [d(\tau)]^{-\frac{1}{p_{+}}}
    E_{\tau}(f)E_{\tau}(g).
\]

We next turn to the loss term. For its $p_{+}$ component, Lemma~\ref{lem:collision-velocity} and Lemma~\ref{lem:collision-frequency} yield
\begin{align*}
    \sum_{|I|\leq 8}
    \norm{\mathcal{D}^{I}Q^{-}(f,g)}_{L_{x,\xi}^{p_{+}}}
    &\lesssim
    \sum_{\substack{|I_1|+|I_2|\leq 8\\ |I_1|\leq 5}}
    \norm{\mathcal{D}^{I_1}f}_{L_{x}^{\infty}L_{\xi}^{p_{+}}}
    \norm{A_{\gamma}[\mathcal{D}^{I_2}g]}_{L_{x}^{p_{+}}L_{\xi}^{\infty}}
    \\
    &\quad
    +
    \sum_{\substack{|I_1|+|I_2|\leq 8\\ |I_2|\leq 2}}
    \norm{\mathcal{D}^{I_1}f}_{L_{x,\xi}^{p_{+}}}
    \norm{A_{\gamma}[\mathcal{D}^{I_2}g]}_{L_{x,\xi}^{\infty}}
    \\
    &\lesssim
    [d(\tau)]^{-\frac{1}{p_{+}}}
    \left\{
    \sum_{\substack{|I_1|+|I_2|\leq 8\\ |I_1|\leq 5,\ |J|\leq 3}}
    \norm{\mathcal{D}^{I_1+J}f}_{L_{x,\xi}^{p_{+}}}
    \left[
    \norm{\mathcal{D}^{I_2}g}_{L_{x}^{p_{+}}L_{\xi}^{p_{-}}}
    +
    \norm{\mathcal{D}^{I_2}g}_{L_{x,\xi}^{p_{+}}}
    \right]
    \right.
    \\
    &\qquad\left.
    +
    \sum_{\substack{|I_1|+|I_2|\leq 8\\ |I_2|\leq 2,\ |J|\leq 3}}
    \norm{\mathcal{D}^{I_1}f}_{L_{x,\xi}^{p_{+}}}
    \left[
    \norm{\mathcal{D}^{I_2+J}g}_{L_{x,\xi}^{p_{-}}}
    +
    \norm{\mathcal{D}^{I_2+J}g}_{L_{x,\xi}^{p_{+}}}
    \right]
    \right\}
    \\
    &\lesssim
    [d(\tau)]^{-\frac{1}{p_{+}}}
    \left\{
    \sum_{\substack{|I_1|+|I_2|\leq 8\\ |I_1|\leq 5,\ |J|\leq 3\\ |K|\leq 3}}
    \norm{\mathcal{D}^{I_1+J}f}_{L_{x,\xi}^{p_{+}}}
    \norm{\mathcal{D}^{I_2+K}g}_{L_{x,\xi}^{p_{-}}}
    +
    E_{\tau}(f)E_{\tau}(g)
    \right\}
    \\
    &\lesssim
    [d(\tau)]^{-\frac{1}{p_{+}}}
    \left\{
    \sum_{\substack{|I_1|\leq 8}}
    \norm{\mathcal{D}^{I_1}f}_{L_{x,\xi}^{p_{+}}}
    \sum_{\substack{|I_2|\leq 11}}
    \norm{\mathcal{D}^{I_2}g}_{L_{x,\xi}^{p_{-}}}
    +
    E_{\tau}(f)E_{\tau}(g)
    \right\}
    \\
    &\lesssim
    [d(\tau)]^{-\frac{1}{p_{+}}}
    E_{\tau}(f)E_{\tau}(g).
\end{align*}

For the $p_{-}$ component, we similarly have
\begin{align*}
    \sum_{|I|\leq 11}
    \norm{\mathcal{D}^{I}Q^{-}(f,g)}_{L_{x,\xi}^{p_{-}}}
    &\lesssim
    \sum_{\substack{|I_1|+|I_2|\leq 11\\ |I_1|\leq 5}}
    \norm{\mathcal{D}^{I_1}f}_{L_{x}^{\infty}L_{\xi}^{p_{c}}}
    \norm{\mathcal{D}^{I_2}g}_{L_{x,\xi}^{p_{-}}}
    \\
    &\quad
    +
    \sum_{\substack{|I_1|+|I_2|\leq 11\\ |I_2|\leq 5}}
    \norm{\mathcal{D}^{I_1}f}_{L_{x,\xi}^{p_{-}}}
    \norm{A_\gamma[\mathcal{D}^{I_2}g]}_{L_{x,\xi}^{\infty}}
    \\
    &\lesssim
    [d(\tau)]^{-\frac{1}{p_{+}}}
    \left\{
    \sum_{\substack{|I_1|+|I_2|\leq 11\\ |I_1|\leq 5,\ |J|\leq 3}}
    \norm{\mathcal{D}^{I_1+J}f}_{L_{x,\xi}^{p_c}}
    \norm{\mathcal{D}^{I_2}g}_{L_{x,\xi}^{p_{-}}}
    \right.
    \\
    &\qquad\left.
    +
    \sum_{\substack{|I_1|+|I_2|\leq 11\\ |I_2|\leq 5,\ |J|\leq 3}}
    \norm{\mathcal{D}^{I_1}f}_{L_{x,\xi}^{p_{-}}}
    \left[
    \norm{\mathcal{D}^{I_2+J}g}_{L_{x,\xi}^{p_{+}}}
    +
    \norm{\mathcal{D}^{I_2+J}g}_{L_{x,\xi}^{p_{-}}}
    \right]
    \right\}
    \\
    &\lesssim
    [d(\tau)]^{-\frac{1}{p_{+}}}
    E_{\tau}(f)E_{\tau}(g).
\end{align*}
Combining the gain and loss estimates proves \eqref{Estimate-Q-X}.

\medskip
\noindent\textbf{Step 2: Homogeneous interactions.}

Since both velocity exponents are controlled through order eleven in $\mathcal H_\tau$, for each
$p\in\{p_{-},p_{+}\}$, Lemma~\ref{lem:collision-velocity} together with interpolation yields
\begin{align*}
    \sum_{|I|\leq 11}
    \norm{\mathcal{D}^{I}Q^{+}(F,G)}_{L_{\xi}^{p}}
    &\lesssim
    \sum_{|I_1|+|I_2|\leq 11}
    \norm{Q^{+}(\mathcal{D}^{I_1}F,\mathcal{D}^{I_2}G)}_{L_{\xi}^{p}}
    \\
    &\lesssim
    \sum_{|I_1|+|I_2|\leq 11}
    \norm{\mathcal{D}^{I_1}F}_{L_{\xi}^{p}}
    \norm{\mathcal{D}^{I_2}G}_{L_{\xi}^{p_{c}}}
    \\
    &\lesssim
    \mathcal{H}_{\tau}(F)\mathcal{H}_{\tau}(G).
\end{align*}
For the loss term, we similarly obtain
\begin{equation}\label{Q-H}
\begin{aligned}
    \sum_{|I|\leq 11}
    \norm{\mathcal{D}^{I}Q^{-}(F,G)}_{L_{\xi}^{p}}
    &\lesssim
    \sum_{|I_1|+|I_2|\leq 11}
    \norm{Q^{-}(\mathcal{D}^{I_1}F,\mathcal{D}^{I_2}G)}_{L_{\xi}^{p}}
    \\
    &\lesssim
    \sum_{|I_1|+|I_2|\leq 11}
    \norm{\mathcal{D}^{I_1}F}_{L_{\xi}^{p}}
    \norm{A_{\gamma}[\mathcal{D}^{I_2}G]}_{L_{\xi}^{\infty}}
    \\
    &\lesssim
    \sum_{|I_1|+|I_2|\leq 11}
    \norm{\mathcal{D}^{I_1}F}_{L_{\xi}^{p}}
    \left(
    \norm{\mathcal{D}^{I_2}G}_{L_{\xi}^{p_{-}}}
    +
    \norm{\mathcal{D}^{I_2}G}_{L_{\xi}^{p_{+}}}
    \right)
    \\
    &\lesssim
    \mathcal{H}_{\tau}(F)\mathcal{H}_{\tau}(G).
\end{aligned}
\end{equation}

Combining the gain and loss estimates gives \eqref{Estimate-Q-X-Hom}. Notice that no spatial Sobolev embedding is needed in this case.

\medskip
\noindent\textbf{Step 3: Mixed interactions.}

Since the homogeneous factor is independent of $x$, the gain estimates together with interpolation immediately give
\begin{align}\label{3-esti}
    E_{\tau}(Q^{+}(f,G))
    &\lesssim
    E_{\tau}(f)\mathcal{H}_{\tau}(G),
    \notag\\
    E_{\tau}(Q^{+}(F,g))
    &\lesssim
    E_{\tau}(g)\mathcal{H}_{\tau}(F),
    \\
    E_{\tau}(Q^{-}(f,G))
    &\lesssim
    E_{\tau}(f)\mathcal{H}_{\tau}(G).
    \notag
\end{align}
It remains only to estimate $Q^{-}(F,g)$. Recalling \eqref{Q-H} with $p=p_{+}$ and arguing in the same way, Lemma~\ref{lem:banach-valued-sobolev} gives
\begin{equation}\label{4-esti}
   \begin{aligned}
    \sum_{|I|\leq 8}
    \norm{\mathcal{D}^{I}Q^{-}(F,g)}_{L_{x,\xi}^{p_{+}}}
    \lesssim&
    \sum_{|I_1|\leq 8,\ |I_2|\leq 8}
    \norm{\mathcal{D}^{I_1}F}_{L_{\xi}^{p_{+}}}
    \left(
    \norm{\mathcal{D}^{I_2}g}_{L_{x}^{p_{+}}L_{\xi}^{p_{-}}}
    +
    \norm{\mathcal{D}^{I_2}g}_{L_{x,\xi}^{p_{+}}}
    \right)\notag\\
    \lesssim&
    \mathcal{H}_{\tau}(F)E_{\tau}(g);\\
    \sum_{|I|\leq 11}
    \norm{\mathcal{D}^{I}Q^{-}(F,g)}_{L_{x,\xi}^{p_{-}}}
    \lesssim&\sum_{|I_1|\leq 11,|I_2|\leq 11}\norm{\mathcal{D}^{I_1}F}_{L_{\xi}^{p_{c}}}\norm{\mathcal{D}^{I_2}g}_{L_{x,\xi}^{p_{-}}}\\
    \lesssim & \mathcal{H}_\tau(F)E_{\tau}(g).
\end{aligned} 
\end{equation}
Therefore, \eqref{3-esti} and \eqref{4-esti} imply
\eqref{Estimate-Q-X-interaction-1} and \eqref{Estimate-Q-X-interaction-2}, respectively.
This completes the proof of Lemma \ref{thm:bilinear-estimation}.
\end{proof}
\begin{remark}
     When $\gamma=0$, the cutoff collision estimate is
$\|Q(u,v)\|_{L^1_\xi}\lesssim\|u\|_{L^1_\xi}\|v\|_{L^1_\xi}$.
Among six distributed derivatives, one input has at most three derivatives. Putting it in $L_x^\infty$ costs three more derivatives and gives $d^{-1}$ by \eqref{Es-Disper}. Homogeneous and mixed interactions follow directly. Thus the six-derivative norm closes at the Maxwell endpoint as well.
\end{remark}

The following lemma is concerned with the estimates of the $Z-$norm of the bilinear operator $Q(\cdot,\cdot)$, where $Z-$norm is given by Definition \ref{defi-shear-hom}.
\begin{lemma}[Collision estimates for homogeneous shear]\label{bi-esti-shear}
Let $-2\leq\gamma\leq 0$, then for all $f,g\in Z^m$ and integer $m\geq 0$:
\begin{equation*}
    \norm{Q(f,g)}_{Z^m} \lesssim \norm{f}_{Z^m} \norm{g}_{Z^m}.
\end{equation*}
\end{lemma}

\begin{proof}
Firstly, we deal with $m=0$. Recall that we denote
$$A_{\gamma}[h](v)
=\displaystyle\int_{\mathbb{R}^3} |v-v_*|^\gamma h(v_*) dv_*.$$
For the loss term,
\[
\norm{Q^-(f,g)}_{L^{3,1}_v}\lesssim \norm{A_\gamma[g]}_{L^\infty_v}\norm{f}_{L^{3,1}_v}.
\]
Moreover, the pre-postcollisional change of variables yields
\[
\begin{aligned}
\|Q^+(f,g)\|_{L^1_v}
&\lesssim \min\bigl\{\|A_\gamma[|g|]\|_{L^\infty_v}\|f\|_{L^1_v},
\|A_\gamma[|f|]\|_{L^\infty_v}\|g\|_{L^1_v}\bigr\}.\\
\|Q(f,g)\|_{L^1_v}
&\lesssim \|A_\gamma[|g|]\|_{L^\infty_v}\|f\|_{L^1_v},
\end{aligned}
\]
For the gain term, a direct change of variables gives
\[
\|Q^+(f,g)\|_{L^\infty_v}\lesssim
\|A_\gamma[|g|]\|_{L^\infty_v}\|f\|_{L^\infty_v}+\|A_\gamma[|f|]\|_{L^\infty_v}\|g\|_{L^\infty_v}.
\]
By positivity, it suffices to consider $f,g\geq0$. Notice that
$\|h\|_{L^{3,1}_v}\lesssim \|h\|_{L^1_v}^{1/3}\|h\|_{L^\infty_v}^{2/3}$,
hence the layer-cake formula and the preceding $L^1_v$ and $L^\infty_v$ estimates yield
\begin{align*}
\|Q^+(f,g)\|_{L^{3,1}_v}
&\leq \int_0^\infty\int_0^\infty\|Q^+(\chi_{\{f>\lambda\}},\chi_{\{g>\mu\}})\|_{L^{3,1}_v}d\mu d\lambda\\
&\lesssim \int_0^\infty\int_0^\infty
\min\Bigl\{|\{f>\lambda\}|\|A_\gamma[\chi_{\{g>\mu\}}]\|_{L^\infty_v},
|\{g>\mu\}|\|A_\gamma[\chi_{\{f>\lambda\}}]\|_{L^\infty_v}\Bigr\}^{1/3}\\
&\qquad\times
\Bigl(\|A_\gamma[\chi_{\{g>\mu\}}]\|_{L^\infty_v}+\|A_\gamma[\chi_{\{f>\lambda\}}]\|_{L^\infty_v}\Bigr)^{2/3}
d\mu d\lambda\\
&\lesssim \int_0^\infty\int_0^\infty
\min\Bigl\{|\{f>\lambda\}|\,|\{g>\mu\}|^{1+\gamma/3},|\{g>\mu\}|\,|\{f>\lambda\}|^{1+\gamma/3}\Bigr\}^{1/3}\\
&\qquad\times
\Bigl(|\{f>\lambda\}|^{1+\gamma/3}+|\{g>\mu\}|^{1+\gamma/3}\Bigr)^{2/3}
d\mu d\lambda\\
&\lesssim \int_0^\infty\int_0^\infty
\Bigl(|\{f>\lambda\}|^{1+\gamma/3}|\{g>\mu\}|^{1/3}+|\{g>\mu\}|^{1+\gamma/3}|\{f>\lambda\}|^{1/3}\Bigr)
d\mu d\lambda\\
&\lesssim\bigl(\|f\|_{L^1_v}+\|f\|_{L^{3,1}_v}\bigr)\|g\|_{L^{3,1}_v}+\|f\|_{L^{3,1}_v}
\bigl(\|g\|_{L^1_v}+\|g\|_{L^{3,1}_v}\bigr).
\end{align*}
Here, the last inequality uses
$\int_0^\infty|\{f>\lambda\}|^{1/3}d\lambda \sim\|f\|_{L^{3,1}_v}$ and $-2\leq \gamma\leq 0$.
Combining the gain and loss estimates with
Lemma~\ref{Endpoint-Loss term} completes the proof for $m=0$.

Secondly, we prove the general $m\in\mathbb{N}^+$. For \(f,g\in Z^m\) and every multi-index \(|\alpha|\leq m\),  we obtain
\[
\|\partial_v^\alpha Q(f,g)\|_Z
\lesssim \sum_{\eta\leq\alpha}
\|\partial_v^\eta f\|_Z \|\partial_v^{\alpha-\eta}g\|_Z.
\]
Summing over \(|\alpha|\leq m\), we complete the proof of
lemma \ref{bi-esti-shear}.
\end{proof}

\section{Dilatation: proof of Theorem~\ref{thm:dilatation}}\label{sec-dial-cylind}
In this section, we shall prove the homogeneous existence and perturbative stability statements together. The two dilatational regimes have the same argument: the trace of the scaled deformation is integrable, and the rescaled collision coefficient belongs to $L^1$ in time. Their different mixed-norm decay rates come from the determinant of the spatial part of the commuting fields.

We first verify the scaling, including the scalar factor that makes the
velocity drift conservative. Recall $a = \text{tr}P/3$, at fixed $(x,v)$, put  $\tau = \log (1+t), \xi = (1+t)^{a}v,$  while
$\nabla_v=(1+t)^a\nabla_\xi$. Therefore $f(t,x,v)=m(t)F(\tau,x,\xi)$
satisfies
\[
 \begin{aligned}
 \frac{1+t}{m}\partial_tf
 &=\partial_\tau F+a\xi\cdot\nabla_\xi F
       +(1+t)\frac{m'}mF,\\
 \frac {1+t}{m}v\cdot\nabla_xf
 &=(1+t)^{1-a}\xi\cdot\nabla_xF,\\
 \frac{1+t}{m}(L(t)v)\cdot\nabla_vf
 &=(P+(1+t)R(t))\xi\cdot\nabla_\xi F.
 \end{aligned}
\]
The choice of $m$ in \eqref{eq:dilatation-scaling} gives
$(1+t) m'/m=-\operatorname{tr}A(\tau)$.
Since $\operatorname{tr}(P-aI)=0$. Consequently, the velocity
terms combine exactly into $-\nabla_\xi\cdot(A\xi F)$.
The collision term has a different scaling:
\[
 Q(f,f)(t,x,v)=m(t)^2(1+t)^{-a(3+\gamma)}Q(F,F)(\tau,x,\xi).
\]
After multiplication by $(1+t)/m$, the coefficient is precisely
$k=M(1+t)^{1-a(3+\gamma)}$. The same computation without spatial transport
gives the homogeneous equation.

Use the scaling \eqref{eq:dilatation-scaling} and the coefficients \eqref{eq:scaled-coefficients}. The system \eqref{eq:scaled-system} reads
\begin{align}
 \mathcal LF=kQ(F,F),\qquad \mathcal L_HG=kQ(G,G),
 \qquad \mathcal T=e^{(1-a)\tau}I. \label{com-eqn}
\end{align}
Thus $\|r(\tau)\|\lesssim e^{-\tau}$, $M$ is bounded above and below by positive constants, and $\operatorname{tr}A=\operatorname{tr}r$. For diagonal dilatation $P=I$, so $a=1$ and $k=Me^{-(2+\gamma)\tau}$. For cylindrical dilatation
\[
 P=\begin{pmatrix}1&0&K\\0&1&0\\0&0&0\end{pmatrix},
 \quad P^2=P,\quad a=\frac23,\quad
 k=Me^{-(1+2\gamma/3)\tau}.
\]
In their respective ranges $-2<\gamma\leq0$ and $-3/2<\gamma\leq0$,
\begin{equation*}
 \int_0^\infty |\operatorname{tr}A(\tau)|\,d\tau<\infty,
 \qquad  \int_0^\infty k(\tau)\,d\tau<\infty.
\end{equation*}

Let $U(\tau,s)$ and $U_H(\tau,s)$ be the linear solution operators for $\mathcal L$ and $\mathcal L_H$. The characteristic map from $s$ to $\tau$ has velocity matrix $\Phi(\tau)\Phi(s)^{-1}$ and determinant
\[
 \frac{\det \Phi(\tau)}{\det \Phi(s)}
 =\exp\left(-\int_s^\tau\operatorname{tr}A(r)\,dr\right).
\]
Solving for the initial point gives the explicit formula
\begin{equation*}
 \begin{aligned}
 (U(\tau,s)h)(x,\xi)
 =\frac{\det \Phi(s)}{\det \Phi(\tau)}
 h\big(&x-(\Psi(\tau)-\Psi(s))\Phi(\tau)^{-1}\xi,
       \Phi(s)\Phi(\tau)^{-1}\xi\big),
 \end{aligned}
\end{equation*}
it follows, for $1\leq p<\infty$, that
\begin{equation}\label{eq:linear-Lp-bound}
 \|U(\tau,s)h\|_{L_{x,\xi}^p}
\lesssim  \|h\|_{L_{x,\xi}^{p}},\qquad  \|U_H(\tau,s)h\|_{L_{\xi}^p}
\lesssim  \|h\|_{L_{\xi}^{p}}.
\end{equation}
And the commuting identities, obviously, imply
\begin{equation}\label{Exchange-for-U-U}
    \mathcal D^I(\tau)U(\tau,s)=U(\tau,s)\mathcal D^I(s),
 \qquad Y_H^\beta(\tau)U_H(\tau,s)=U_H(\tau,s)Y_H^\beta(s).
\end{equation}

Consequently, the linear operators are uniformly bounded in the adapted energies. 

We next give the following elementary lemma concerning the determinant of the weight for vector fields
\begin{lemma}\label{lem:matrix-growth}
For both dilatational regimes, there exist positive constants $C$ and $c$ such that
$$0<c\leq\det \Phi(\tau)\leq C.$$  
Moreover, there exists $\tau_0>0$ such that, for all $\tau\geq\tau_0$,
\[
 \det \Psi(\tau)\sim
 \begin{cases}
 (1+\tau)^3,&P=I,\\
 e^\tau(1+\tau)^2,&P=P_K.
 \end{cases}
\]
\end{lemma}

\begin{proof}
We consider the two cases separately.

\medskip
\noindent
\textbf{Case 1: Homogeneous diagonal.}
In this case, $\mathcal{T}(\tau)=I$ and
\[
    A(\tau)=r(\tau), \qquad \norm{r(\tau)}\lesssim e^{-\tau}.
\]
By \eqref{B-C-defi},
\[
    \Phi(\tau)=\Phi(0)-\int_0^\tau r(s)\Phi(s)\,ds.
\]
Gronwall's inequality gives
\[
    \norm{\Phi(\tau)}
    \leq \norm{\Phi(0)}
    \exp\left\{\int_0^\tau\norm{r(s)}\,ds\right\}
    \lesssim 1.
\]
Consequently,
\[
    \int_s^t\norm{\Phi'(\tau)}\,d\tau
    \lesssim\int_s^t\norm{r(\tau)}\,d\tau,
\]
and hence $\Phi(\tau)$ converges as $\tau\to\infty$. We denote
 $ \displaystyle
    \Phi(\infty)=\lim_{\tau\to\infty}\Phi(\tau).
$

By Liouville's formula, it follws that 
\begin{align*}
    \det\Phi(\tau)
    &=\det\Phi(0)
      \exp\left\{-\int_0^\tau\operatorname{tr}r(s)\,ds\right\}\neq0,\\
    \det\Phi(\infty)
    &=\det\Phi(0)
      \exp\left\{-\int_0^\infty\operatorname{tr}r(s)\,ds\right\}\sim1.
\end{align*}
Thus,
\(
    \det\Phi(\tau)\sim1.
\)
Again by \eqref{B-C-defi},
\[
    \frac{\Psi(\tau)}{\tau}-\Phi(\infty)
    =\frac{\Psi(0)}{\tau}
    +\frac1\tau\int_0^\tau
      \bigl(\Phi(s)-\Phi(\infty)\bigr)\,ds.
\]
Since
\[
    \norm{\Phi(\tau)-\Phi(\infty)}
    =\norm{\int_\tau^\infty r(s)\Phi(s)\,ds}
    \lesssim e^{-\tau},
\]
we obtain
\[
    \norm{\frac{\Psi(\tau)}{\tau}-\Phi(\infty)}
    \lesssim\frac1\tau.
\]
Writing
\[
    \Psi(\tau)=\tau\bigl(\Phi(\infty)+\beta(\tau)\bigr),
    \qquad
    \norm{\beta(\tau)}\lesssim\tau^{-1}\to0,
\]
we have
\[
    \det\Psi(\tau)
    =\tau^3\det\bigl(\Phi(\infty)+\beta(\tau)\bigr).
\]
By the determinant perturbation estimate
\[
    \left|
    \det\bigl(\Phi(\infty)+\beta(\tau)\bigr)
    -\det\Phi(\infty)
    \right|
    \leq
    3\bigl(\norm{\Phi(\infty)}+\norm{\beta(\tau)}\bigr)^2
    \norm{\beta(\tau)},
\]
one further obtains
\[
    \lim_{\tau\to\infty}
    \frac{\det\Psi(\tau)}{\tau^3}
    =\det\Phi(\infty),
\]
and hence
\[
    \det\Psi(\tau)
    \sim\tau^3\det\Phi(\infty)
    \gtrsim d(\tau).
\]

\medskip
\noindent
\textbf{Case 2: Cylindrical dilatation with shear.}
In this case,
\[
    A(\tau)
    =P-\frac23I+r(\tau),
    \qquad
    \norm{r(\tau)}\lesssim e^{-\tau}.
\]

We first estimate $\det\Phi(\tau)$. Since
$\Phi'=-A\Phi$, it follows
\[
    \det\Phi(\tau)
    =
    \exp\left\{
        -\int_0^\tau\operatorname{tr}A(s)\,ds
    \right\}\det\Phi(0)
    =
    \exp\left\{
        -\int_0^\tau\operatorname{tr}A(s)\,ds
    \right\}.
\]
Thus
\[
    \det\Phi(\tau)\sim1.
\]

We next estimate $\det\Psi(\tau)$. Since $\Phi(\tau)$ is invertible,
set
\[
    J(\tau)
    =e^{-\tau/3}\Psi(\tau)\Phi(\tau)^{-1},
    \qquad
    \det J(\tau)
    =\frac{e^{-\tau}\det\Psi(\tau)}
           {\det\Phi(\tau)}.
\]
Since
\[
    (\Phi^{-1})'
    =-\Phi^{-1}\Phi'\Phi^{-1}
    =
    \Phi^{-1}
    \left(P-\frac23I+r(\tau)\right),
\]
we obtain
\begin{equation}\label{eq-for-H(tau)}
    J'=-J(I-P)+I+Jr(\tau).
\end{equation}
Using $\norm{r(\tau)}\lesssim e^{-\tau}$ and Gronwall's inequality,
\[
    \norm{J(\tau)}
    \lesssim
    (1+\tau)
    +\int_0^\tau\norm{J(s)}e^{-s}\,ds
    \lesssim1+\tau.
\]

Multiplying \eqref{eq-for-H(tau)} by $P$ and using
$P(I-P)=0$, we obtain
\[
    (JP)'=P+Jr(\tau)P,
\]
and hence
\begin{equation}\label{HP-esti}
    JP=(1+\tau)P+W_1(\tau),
    \qquad
    W_1(\tau)
    =\int_0^\tau J(s)r(s)P\,ds,
    \qquad
    \norm{W_1(\tau)}\lesssim1.
\end{equation}
Similarly, multiplying \eqref{eq-for-H(tau)} by $I-P$ gives
\[
    (J(I-P))'
    =I-P-J(I-P)+Jr(\tau)(I-P),
\]
and therefore
\begin{equation}\label{est I - P}
    J(\tau)(I-P)
    =I-P+W_2(\tau),
    \qquad
    \norm{W_2(\tau)}
    \lesssim e^{-\tau}\tau(1+\tau).
\end{equation}
Combining \eqref{HP-esti} and \eqref{est I - P}, we obtain
\[
    J\left(\frac{P}{1+\tau}+I-P\right)
    =
    I+\frac{W_1(\tau)}{1+\tau}+W_2(\tau).
\]
Since
\[
    \frac{W_1(\tau)}{1+\tau}+W_2(\tau)\to0,
\]
it follows that
\[
    \det\left(
        I+\frac{W_1(\tau)}{1+\tau}+W_2(\tau)
    \right)\to1.
\]
Thus, for all sufficiently large $\tau$,
\[
    \det J(\tau)
    \det\left(\frac{P}{1+\tau}+I-P\right)
    \geq\frac12.
\]
Since
\[
    \frac{P}{1+\tau}+I-P
    =
    \begin{pmatrix}
        (1+\tau)^{-1}&0&K(1+\tau)^{-1}-K\\
        0&(1+\tau)^{-1}&0\\
        0&0&1
    \end{pmatrix},
\]
we have
\[
    \det\left(\frac{P}{1+\tau}+I-P\right)
    =(1+\tau)^{-2},
\]
and hence
\[
    \det J(\tau)\geq\frac{(1+\tau)^2}{2}.
\]
Finally, using $\det\Phi(\tau)\sim1$,
\[
    \det\Psi(\tau)
    =e^\tau\det J(\tau)\det\Phi(\tau)
    \gtrsim e^\tau\tau^2
    \sim d(\tau).
\]
This completes the proof of Lemma \ref{lem:matrix-growth}.
\end{proof}


We are now ready to complete the proof of Theorem~\ref{thm:dilatation}. We focus on the case $\gamma<0$. The Maxwell case $\gamma=0$ can be treated by a similar but simpler argument, relying on sixth-order $L^1$ estimates, and we therefore omit the details.

\begin{proof}[The proof of Theorem~\ref{thm:dilatation}]
We first solve the homogeneous equation by its Duhamel formula
\[
 G(\tau)=U_H(\tau,0)G_0+
       \int_0^\tau U_H(\tau,s)k(s)Q(G,G)(s) ds.
\]
Applying vector fields $Y_H^{\beta}$, and \eqref{Exchange-for-U-U} yields
\begin{equation}\label{vec-mild-formu}
       Y_H^{\beta}G(\tau)=U_H(\tau,0)[Y_H^{\beta}(0)G_0]+
       \int_0^\tau U_H(\tau,s)[k(s)Y_H^{\beta}(s)Q(G,G)(s)] ds.
\end{equation}
Recall \eqref{eq:collision-time-integrability}, then bilinear estimate~\eqref{Estimate-Q-X-Hom} and \eqref{eq:linear-Lp-bound} imply, uniformly in $T$,
    \begin{align*}
   \sum_{|\beta|\leq 11} \norm{Y_H^{\beta}G(\tau)}_{L_{\xi}^p}\lesssim&\sum_{|\beta|\leq 11} \norm{Y_H^{\beta}(0)G_0}_{L_{\xi}^{p}} + \sum_{|\beta|\leq 11} \int_{0}^{T}k(s)\norm{Y_H^{\beta}Q(G,G)(s)}_{L_{\xi}^{p}}ds\\
    \lesssim&\sum_{|\beta|\leq 11} \norm{Y_H^{\beta}(0)G_0}_{L_{\xi}^{p}}+\int_{0}^{T}k(s)[\mathcal{H}_{s}(G)]^2ds\\\lesssim&\mathcal{H}_{0}(G_0) + [\mathcal{H}_{T}(G)]^2,
    \end{align*}
Summing over $p=p_{+}$ and $p_{-}$, therefore
\begin{equation}\label{Hom-priori}
    \mathcal{H}_{T}(G)\lesssim \mathcal{H}_{0}(G_0) + [\mathcal{H}_{T}(G)]^2.
\end{equation}
Having established the a priori estimates \eqref{Hom-priori} above, the same estimates, together with a standard contraction argument, yield a local mild solution. For sufficiently small initial data, the a priori bound and a standard bootstrap argument allow us to extend this local solution globally in time, with
\begin{equation}\label{smallness}
     \mathcal H_\infty(G)\lesssim \mathcal H_0(G_0).
\end{equation}

Applying the same estimates to the difference of two solutions also shows that the solution map
   $ G_0\longmapsto G$
is Lipschitz continuous. Finally, uniqueness follows from the difference estimate
\[
    \mathcal H_T(f-g)
    \lesssim
    \bigl(\mathcal H_T(f)+\mathcal H_T(g)\bigr)
    \mathcal H_T(f-g),
\]
with  $\mathcal H_T(f)+\mathcal H_T(g)\lesssim \mathcal{H}_0(f_0)\ll 1$, uniformly in $T$.

We now turn to the corresponding inhomogeneous case. Denote $F=G+H$. The perturbation satisfies
\begin{equation*}
\begin{aligned}
     \mathcal LH&=k\mathcal N(H,G),\
 \mathcal N(H,G)=Q(H,G)+Q(G,H)+Q(H,H),\ H(0)=H_0.
\end{aligned}
\end{equation*}
Arguing  as in  \eqref{vec-mild-formu}, bilinear estimate \eqref{Estimate-Q-X}, \eqref{Estimate-Q-X-interaction-1}, \eqref{Estimate-Q-X-interaction-2}  imply
\begin{equation}\label{nonhom-priori-esti-pre}
    \begin{aligned}
   \sum_{|I|\leq N_p} \norm{\mathcal{D}^{I}H}_{L_{x,\xi}^p}\lesssim&
   \sum_{|I|\leq N_p}\norm{\mathcal{D}^{I}(0)H_0}_{L_{x,\xi}^p} + \int_{0}^{T}k(s) \sum_{|I|\leq N_p}\bigg[\norm{\mathcal{D}^{I}Q(H,H)}_{L_{x,\xi}^p}\\& + \norm{\mathcal{D}^{I}Q(G,H)}_{L_{x,\xi}^p}+\norm{\mathcal{D}^{I}Q(H,G)}_{L_{x,\xi}^p} \bigg]ds\\
   \lesssim&E_0(H_0) + \int_{0}^{T}k(s)\left([E_{s}(H)]^2 + E_{s}(H)\mathcal{H}_s (G) \right)ds\\
   \lesssim&E_0(H_0) + [E_{T}(H)]^2 + E_{T}(H)\mathcal{H}_T (G).
\end{aligned}
\end{equation}
Summing over $p = p_{+}$ and $p=p_{-}$, \eqref{nonhom-priori-esti-pre} yields
\begin{equation}\label{nonhom-priori}
  E_T(H) \lesssim  E_0(H_0) + [E_{T}(H)]^2 + E_{T}(H)\mathcal{H}_T (G).
\end{equation}

For sufficiently small initial data $\mathcal H_0(G_0)+E_0(H_0)\ll 1$ , \eqref{smallness} yields $\mathcal{H}_{T}(G)\ll 1$, therefore \eqref{nonhom-priori} implies
\begin{equation}\label{priori-2}
 E_T(H) \lesssim  E_0(H_0) + [E_{T}(H)]^2.
\end{equation}
Together with
a standard contraction argument and a priori estimate \eqref{priori-2}, the standard bootstrap argument  allow us to extend this
local solution globally in time, with
\begin{equation}\label{eq:H-uniform-bound}
    E_\infty(H) \lesssim E_0(H_0).
\end{equation}
  Moreover, the solution map $H_0 \mapsto H$ is Lipschitz continuous, and the uniqueness follows the difference estimate
$$ E_{T}(H_1 - H_2)\lesssim E_0(H_{1,0} - H_{2,0}) + (E_{T}(H_1)+E_T(H_2))E_{T}(H_1 - H_2) + \mathcal{H}_T(G)E_{T}(H_1 - H_2).$$
Notice that the solution belongs to the affine class $G+X$: a nonzero homogeneous background does not belong to $L^p_{x,\xi}(\mathbb R^6)$.

For $G_0\geq0$, the gain--loss formula along characteristics gives $G\geq0$.
If also $G_0\in L^1_\xi$, the collision-frequency bound propagates $L^1$
integrability on finite time intervals. Since both the linear flow and
collisions conserve mass, we obtain
$$
\int_{\mathbb R^3}G(\tau,\xi)\,d\xi
=\int_{\mathbb R^3}G_0(\xi)\,d\xi,
\qquad \tau\geq0.
$$

As $G\geq0$ and $F_0=G_0+H_0\geq0$, the reconstructed solution
$F=G+H$ remains nonnegative throughout its interval of existence,
although the perturbation $H$ may change sign.
Indeed, the bound $F^{-}\leq|H|$ and the collision estimates yield
a Gronwall inequality for the negative part $F^{-}$, which vanishes
initially and therefore remains zero.
We omit the details of this standard positivity argument.

Finally, applying \eqref{Es-Disper} to the differentiated perturbation and using \eqref{eq:H-uniform-bound} gives, for $\gamma<0$,
\begin{align*}
 \sum_{|I|\leq5}\|\mathcal D^I H(\tau)\|_{L_x^\infty L_\xi^{p_+}}
 &\lesssim [d(\tau)]^{-1/p_+}E_0(H_0),\\
 \sum_{|I|\leq8}\|\mathcal D^I H(\tau)\|_{L_x^\infty L_\xi^{p_-}}
 &\lesssim [d(\tau)]^{-1/p_-}E_0(H_0).
\end{align*}
For $\gamma=0$, the corresponding estimate is
\[
 \sum_{|I|\leq3}\|\mathcal D^I H(\tau)\|_{L_x^\infty L_\xi^1}
 \lesssim [d(\tau)]^{-1}E_0(H_0).
\]
These estimates prove the stated convergence. In particular, the cylindrical bound contains the factor $e^{-\tau/p}(1+\tau)^{-2/p}$, and hence implies the asserted exponential bound in scaled time. This ends the proof of Theorem~\ref{thm:dilatation}.
\end{proof}
Finally, we conclude with the following remark on the decay rates of the solutions in the original variables.
\begin{remark}
The zeroth-order decay can also be written in the original variables.
Indeed, $f(t,x,v)-g(t,v)=m(t)H(\log(1+t),x,(1+t)^av)$, so a change of
velocity variable gives, for $p$ as in Theorem~\ref{thm:dilatation},
\[
 \|f(t)-g(t)\|_{L_x^\infty L_v^p}
 =m(t)(1+t)^{-3a/p}\|H(\log(1+t))\|_{L_x^\infty L_\xi^p}.
\]
Since $m$ is uniformly bounded, the theorem implies
\[
 \|f(t)-g(t)\|_{L_x^\infty L_v^p}
 \lesssim E_0(H_0)(1+t)^{-3/p}
 \begin{cases}
 [1+\log(1+t)]^{-3/p},&P=I,\\
 [1+\log(1+t)]^{-2/p},&P=P_K.
 \end{cases}
\]
Thus the different rates in scaled time are compatible with a common
algebraic factor in physical time; the velocity change contributes
to this factor as well. These are mixed-norm bounds for the same
spatial coordinate used in the transport--deformation model.
The higher derivative estimates are most naturally kept in their
commuting-field form, because conversion to ordinary derivatives
introduces time-dependent matrices. Likewise, the uniform bound on
$E_\tau(H)$ is consistent with convergence of the lower mixed norms:
the latter expresses the spreading measured by $\det \Psi$.
\end{remark}

\section{Simple shear: proof of Theorem~\ref{thm:shear}}\label{sec-5}
In this section, we first prove the homogeneous assertion in
Theorem~\ref{thm:shear}\textup{(\ref{part:shear-homogeneous})}.
Since spatial dispersion is absent, we combine collision dissipation
of Fisher information with its growth under shear. We justify this
estimate by a positive splitting scheme and then pass from bounded
collision kernels to soft potentials. In particular, the Lorentz norm was introduced to handle the case of $\ga=-2$ which is not included in the previous dilatation cases. We then prove the inhomogeneous
assertion in Theorem~\ref{thm:shear}\textup{(\ref{part:shear-inhomogeneous})},
where transport dispersion makes the quadratic collision term
integrable in time for small data.

\subsection{Collision dissipation and the shear contribution}
We first establish the two ingredients of the homogeneous Fisher
estimate. In this subsection we denote $b(c)=\frac12(B_0(c)+B_0(-c)).$ And for a nonnegative pair density $R=R(v,v_*)$, write
\[
 \mathcal I_2(R)
 =4\iint_{\mathbb R^6}
 \bigl(|\nabla_v\sqrt R|^2+|\nabla_{v_*}\sqrt R|^2\bigr)\,dv\,dv_*.
\]

\begin{lemma}[Fisher contraction for bounded collision kernels]
\label{prop:collision-fisher}
Let $-2\leq\gamma\leq0$, $\varepsilon>0$, and
$K_\varepsilon(\rho)=(\rho^2+\varepsilon^2)^{\gamma/2}$.
Define
\[
 (\mathcal L_\varepsilon R)(v,v_*)
 =\int_{\mathbb S^2}K_\varepsilon(|v-v_*|)b(e\cdot\omega)
 [R(v',v_*')-R(v,v_*)]\,d\omega,
 \qquad e=\frac{v-v_*}{|v-v_*|}.
\]
Then $\exp(s\mathcal L_\varepsilon)$ is a positive,
mass-preserving contraction on $L^1_{v,v_*}(\mathbb R^6)$. For every symmetric $R\geq 0$ with finite mass and finite $\mathcal{I}_2(R)$,
\begin{equation}\label{eq:pair-fisher-contraction}
 \mathcal I_2(\exp(s\mathcal L_\varepsilon)R)\leq\mathcal I_2(R),
 \qquad s\geq0.
\end{equation}
\end{lemma}

\begin{proof}
The cutoff assumption and boundedness of $K_\varepsilon$ give a
bounded conservative jump generator on $L^1_{v,v_*}(\mathbb R^6)$. Its
semigroup preserves positivity, mass, and exchange symmetry.

Next, we give the Fisher calculation to specify precisely the use of
\cite{ISV}.
Set $v=z+r\sigma$, $v_*=z-r\sigma$, and
$V(z,r,\sigma)=R(z+r\sigma,z-r\sigma)$. Then
$V(z,r,\sigma)=V(z,r,-\sigma)$ and
\[
 \mathcal I_2(R)=4(J_z+J_r+J_\sigma),
\]
where
\[
 \begin{aligned}
 J_z&=\iiint\frac{|\nabla_zV|^2}{V}r^2\,d\sigma\,dr\,dz,
 &J_r&=\iiint\frac{|\partial_rV|^2}{V}r^2\,d\sigma\,dr\,dz,\\
 J_\sigma&=\iiint\frac{|\nabla_\sigma V|^2}{V}\,d\sigma\,dr\,dz.
 \end{aligned}
\]
In these coordinates the generator is $a_\varepsilon(r)\mathcal B$,
with $a_\varepsilon(r)=K_\varepsilon(2r)$ and
\[
 \mathcal B u(\sigma)
 =\int_{\mathbb S^2}\widetilde b(\sigma\cdot\sigma')
 [u(\sigma')-u(\sigma)]\,d\sigma',
 \qquad
 \widetilde b(y)
 =\frac{b(\sqrt{(1-y)/2})}{2\sqrt{(1-y)/2}}.
\]
The formula follows from
$\sigma'=\sigma-2(\sigma\cdot\omega)\omega$ and the evenness of $b$.
Moreover, $\widetilde b$ is integrable on the sphere.

For a smooth positive $V$, write
$d\mu=r^2\widetilde b(\sigma\cdot\sigma')\,d\sigma'\,d\sigma\,dr\,dz$,
$U=V(z,r,\sigma)$, and $U'=V(z,r,\sigma')$.
Denote $a'_{\varepsilon}=\frac{d}{dr} a_{\varepsilon}(r)$, then
differentiate the three terms as follow:
\begin{enumerate}
    \item For the parallel term, a direct calculation gives
    \[
    DJ_z(V)[a_\varepsilon\mathcal BV]=-\frac12 \int a_{\varepsilon}
    (U+U')|\nabla_z \log U'-\nabla_z\log U|^2 d\mu \leq 0.
    \]
    \item For the radial term, a direct calculation gives
    \begin{align*}
     &DJ_r(V)[a_\varepsilon\mathcal BV]\\
     &=\int \left(-\frac {a_{\varepsilon}}{2}  (U+U') (\partial_r\log U'-\partial_r\log U)^2 -a_{\varepsilon}'(U'-U)
    (\partial_r\log U'-\partial_r\log U)\right) d\mu \\
    &\leq \frac12\int\frac{|a_\varepsilon'|^2}{a_\varepsilon}
          \frac{(U'-U)^2}{U'+U} d\mu.
    \end{align*}

    \item For the spherical term, the spherical Fisher inequality of
\cite[Proposition~7.11]{ISV}, applied to the even function $u$, gives
\[
 DJ_\sigma(V)[a_\varepsilon\mathcal BV]
 \leq-2\Lambda_{\widetilde b}\int
 \frac{a_\varepsilon}{r^2}\frac{(U'-U)^2}{U'+U}\,d\mu,
 \qquad \Lambda_{\widetilde b}\geq1,
\]
\end{enumerate}
where $\Lambda_{\widetilde b}$ is defined in \cite{ISV}.
Combining (1) (2) (3), we obtain
\[
 \frac{d}{ds}\mathcal I_2(\exp(s\mathcal L_\varepsilon)R)
 \leq4\int\left(\frac{|a_\varepsilon'|^2}{2a_\varepsilon}
 -\frac{2\Lambda_{\widetilde b}a_\varepsilon}{r^2}\right)
 \frac{(U'-U)^2}{U'+U}\,d\mu\leq0,
\]
since
\[
 \frac{r|a_\varepsilon'(r)|}{a_\varepsilon(r)}
 =|\gamma|\frac{4r^2}{4r^2+\varepsilon^2}
 \leq|\gamma|\leq 2\leq2\sqrt{\Lambda_{\widetilde b}}.
\]

To extend the Fisher contraction estimate from smooth data to general pair densities with finite Fisher information as allowed in the lemma,  approximate $\widetilde b$ in angular
$L^1_{\sigma'}$ by nonnegative smooth zonal kernels $\widetilde b_n$. Let $\mathcal L_{\varepsilon,n}$ and $P_{\varepsilon,n}(s)=e^{s\mathcal L_{\varepsilon,n}}$
denote the corresponding generators and semigroups. Then $\mathcal L_{\varepsilon,n} \to \mathcal L_\varepsilon$ in $L^1_{v,v_*}$ operator norm,
and $P_{\varepsilon,n}(s) \to \exp(s\mathcal L_\varepsilon)$ converge strongly.

Choose nonnegative smooth symmetric $W_n$ such that
\[
 W_n=0\quad\text{if }|v-v_*|\leq\eta_n,\qquad
 \eta_n\downarrow0,\qquad
 W_n\longrightarrow\sqrt R\quad\text{in }H^1(\mathbb R^6).
\]
Such a sequence is obtained by cutting off $\sqrt R$ near $v=v_*$, with the cutoff error tending to zero in $H^1$ by the three-dimensional Hardy inequality in $v-v_*$, followed by symmetric compact truncation and nonnegative mollification.
Set
\[
 R_n=W_n^2+c_ne^{-|v|^2-|v_*|^2},\qquad c_n>0,\quad c_n\longrightarrow0.
\]
Then $R_n>0$, $R_n(v,v_*)=R_n(v_*,v)$, and
\[
 \|R_n-R\|_{L^1_{v,v_*}}
 \leq\|W_n-\sqrt R\|_{L^2_{v,v_*}}(\|W_n\|_{L^2_{v,v_*}}+\|\sqrt R\|_{L^2_{v,v_*}})+c_n\|e^{-|v|^2-|v_*|^2}\|_{L^1_{v,v_*}}
 \longrightarrow0.
\]
By subadditivity of Fisher information,
\[
\mathcal I_2(R_n)
\leq4\|\nabla_{v,v_*}W_n\|_2^2+c_n\mathcal I_2(e^{-|v|^2-|v_*|^2})
\longrightarrow\mathcal I_2(R).
\]
Also,
$\|\sqrt{R_n}-\sqrt R\|_{L^2}^2\leq\|R_n-R\|_{L^1}\to0$.
Choose a subsequence realizing $\liminf_n\mathcal I_2(R_n)$.
The preceding Fisher upper bound makes its square roots bounded
in $H^1$; after a further subsequence they converge weakly in
$H^1$ to $\sqrt R$, as identified by the strong $L^2$ limit.
Hence
\[
\mathcal I_2(R)
\leq\liminf_n\mathcal I_2(R_n)
\leq\limsup_n\mathcal I_2(R_n)
\leq\mathcal I_2(R).
\]
Therefore $\mathcal I_2(R_n)\to\mathcal I_2(R)$.
Applying the differentiating calculation above to $\widetilde b_n$, we have
 \[
 \frac{d}{ds}\mathcal I_2(P_{\varepsilon,n}(s)R_n)\leq0, \qquad
 \mathcal I_2(P_{\varepsilon,n}(s)R_n)-\mathcal I_2(R_n)
=\int_0^s\frac{d}{dt}\mathcal I_2(P_{\varepsilon,n}(t)R_n) dt\leq0,
 \]
 which implies,
for every fixed $s\geq0$,
\[
 \|P_{\varepsilon,n}(s)R_n-\exp(s\mathcal L_\varepsilon)R\|_{L^1_{v,v_*}}
 \leq\|R_n-R\|_{L^1_{v,v_*}}+s\|\mathcal L_{\varepsilon,n}
 -\mathcal L_\varepsilon\|\|R\|_{L^1_{v,v_*}} \longrightarrow0,
\]
and then, using the weak convergence argument as above,
\[
 \mathcal I_2(\exp(s\mathcal L_\varepsilon)R)
 \leq\liminf_n\mathcal I_2(P_{\varepsilon,n}(s)R_n)
 \leq\lim_n\mathcal I_2(R_n)=\mathcal I_2(R).
\]
This proves \eqref{eq:pair-fisher-contraction}.
\end{proof}

With the lemma above in hand, we can get a priori estimate on the existence interval.
\begin{proposition}[Fisher bound for the full shear equation]
\label{prop:Fisher}
Let $-2\leq\gamma\leq0$, $g_0\geq0$, $M_0=\norm{g_0}_{L^1_v}<\infty$, and
$\mathcal I(g_0)<\infty$. If $g$ is a mild solution of \eqref{sim-shear-hom} with initial value $g_0$ in $C([0,T];Z)$, that is:
\[
g(t)=S(t)\left[g_0 + \int_0^t S(-s)Q(g(s),g(s)) ds\right].
\]
Then $g\geq 0,\ \norm{g}_{L_v^1}=M_0$ and
\begin{equation}\label{eq:full-Fisher}
 \mathcal I(g(t))\leq e^t \mathcal{I}(g_0),\qquad 0\leq t\leq T.
\end{equation}
\end{proposition}

\begin{proof}
Without loss of generality, suppose $M_0>0$.

\textbf{Step 1: Fisher information of the simple shear term}

For $u=\sqrt{h}\in H^1$, recall the shear group in
Definition~\ref{defi-shear-hom}:
$\sqrt{S(t)h}(v_1,v_2,v_3)=u(v_1+tv_2,v_2,v_3)$.
Then the change of variables $y=(v_1+tv_2,v_2,v_3)$ gives
\[
\begin{aligned}
 \mathcal I(S(t)h)
 &=4\int_{\R^3}
 \left[|\partial_1u|^2+
 |t\partial_1u+\partial_2u|^2+|\partial_3u|^2\right]dy,\\
 \frac{d}{dt}\mathcal I(S(t)h)
 &=8\int\partial_1u(t\partial_1u+\partial_2u)dy
 \leq\mathcal I(S(t)h).
\end{aligned}
\]
Gronwall's inequality therefore yields
\[
\mathcal I(S(t)h)\leq e^t\mathcal I(h),\qquad t\geq0.
\]

\textbf{Step 2: A splitting argument}

Let $Q_\varepsilon$ be \eqref{eq:collision-definition} with $|v-v_*|^\gamma$ replaced by $K_\varepsilon(|v-v_*|)$ defined in lemma \ref{prop:collision-fisher}, and use the semigroup $\exp(s\mathcal L_\varepsilon)$ from
Lemma~\ref{prop:collision-fisher}. For fixed $\varepsilon>0$ and step size $\delta\geq0$, we suppress their dependence in $\mathcal C$. For $f\geq0$ of mass $M_0$, set
\begin{equation*}
\mathcal C[f]
=\frac{1}{M_0}\pi \exp(M_0\delta\mathcal L_\varepsilon)(f\otimes f),
\quad
(\pi R)(v)=\int R(v,v_*)\,dv_*.
\end{equation*}
Notice that for symmetric finite-Fisher $R\geq0$, we have $\mathcal I_2(f\otimes f)=2M_0\mathcal I(f)$, while Cauchy--Schwarz gives
$2\mathcal I(\pi R)\leq\mathcal I_2(R)$. Thus
\begin{equation}\label{eq: approximation bound}
    \mathcal I(\mathcal C[f])
 \leq\frac1{2M_0}
 \mathcal I_2\big(\exp(M_0\delta\mathcal L_\varepsilon)(f\otimes f)\big)
 \leq\frac1{2M_0}\mathcal I_2(f\otimes f)
 =\mathcal I(f).
\end{equation}
Consequently, starting from the initial datum $g_0$, define the splitting sequence by
$g_{j+1}=S(\delta)\mathcal C[g_j]$, $j\geq0$.
then by Step 1 and \eqref{eq: approximation bound}:
\begin{equation}\label{eq:split-Fisher}
 g_j\geq0,\qquad\norm{g_j}_{L^1_v}=M_0,\qquad
 \mathcal I(g_j)\leq e^{j\delta}\mathcal I(g_0).
\end{equation}

Moreover, $\pi\mathcal L_\varepsilon(g\otimes g)
=Q_\varepsilon(g,g)$ and the expansion
$\exp(M_0\delta\mathcal L_\varepsilon)=\mathrm{Id}+M_0\delta\mathcal L_\varepsilon+O(\delta^2)$ imply
\begin{equation}\label{eq:collision-step-stability}
    \norm{\mathcal C[g]-g-\delta Q_\varepsilon(g,g)}_{L^1_v}
 \leq C_{\varepsilon}\delta^2.
\end{equation}
Finally, \eqref{eq:split-Fisher} and Lemma~\ref{Endpoint-Loss term} yield
\begin{equation}\label{eq:split-frequency}
 \sup_{j\delta\leq T}
\left(\norm{g_j}_Z+\norm{A_\gamma[g_j]}_{L^\infty_v}\right)
\lesssim M_0+e^T\mathcal I(g_0),
\end{equation}
uniformly in $\varepsilon$ and $\delta$.

\textbf{Step 3: Proof of the proposition}

Set
$l_j=\|g_j-g(j\delta)\|_{L^1_v}$.
Then \eqref{eq:split-frequency} gives
\[
\norm{Q_\varepsilon(g_j,g_j)-Q_\varepsilon(g(j\delta),g(j\delta))}_{L^1_v}\\
\lesssim\left(\norm{A_\gamma[g_j]}_{L^\infty_v}+\norm{A_\gamma[|g(j\delta)|]}_{L^\infty_v}\right)l_j
\leq C_Tl_j,
\]
where the inequality is independent of $\varepsilon,\delta$.
Thus, by \eqref{eq:collision-step-stability}, whenever $(j+1)\delta\leq T$,
\begin{align*}
l_{j+1}&=\norm{S(\delta)\mathcal C[g_j]-
S(\delta) g(j\delta)-\int_0^\delta S(\delta-s) Q(g(j\delta+s),g(j\delta+s))ds}_{L^1_v} \\
&\leq
\left\|S(\delta)\left\{g_j-g(j\delta)+\delta\left[Q_\varepsilon(g_j,g_j)
 -Q_\varepsilon(g(j\delta),g(j\delta))\right]\right\}\right\|_{L^1_v}\\
&\quad+
\left\|S(\delta)\left\{\mathcal C[g_j]-g_j
 -\delta Q_\varepsilon(g_j,g_j)\right\}\right\|_{L^1_v}\\
&\quad+\delta\norm{S(\delta)\left[Q_\varepsilon(g(j\delta),g(j\delta))
 -Q(g(j\delta),g(j\delta))\right]}_{L^1_v}\\
&\quad+\int_0^\delta
\norm{S(\delta)Q(g(j\delta),g(j\delta))-S(\delta-s)Q(g(j\delta+s),g(j\delta+s))}_{L^1_v}\,ds\\
&\leq (1+C_T\delta)l_j+C_\varepsilon\delta^2
 +\delta\eta_1(\varepsilon)+\delta\eta_2(\delta).
\end{align*}
Here,
\[
\eta_1(\varepsilon):=\sup_{t\leq T}\norm{Q_\varepsilon(g(t),g(t))-Q(g(t),g(t))}_{L^1_v}.
\]
and
\[
\eta_2(\delta):=\sup_{\substack{0\leq t\leq T-\delta\\0\leq s\leq\delta}}
\norm{Q(g(t),g(t))-S(-s)Q(g(t+s),g(t+s))}_{L^1_v}
\]
On one hand, strong continuity of $S$ implies $\eta_2(\delta)\to 0$ as $\delta\to 0$; on the other hand, we claim that $\eta_1(\varepsilon)\to 0$ as $\varepsilon\to 0$: Indeed,
since $g\in C([0,T];Z)$ and $L^{3,1}_v\hookrightarrow L^3_v$,
$\sup_{t\leq T}\norm{g(t)}_{L^2}^2
\leq \sup_{s\leq T}\|g(s)\|_{L^1_v}^{1/2}\sup_{t\leq T}\norm{g(t)}_{L^3_v}^{3/2}<\infty,$ which implies
\[
\begin{aligned}
\eta_1(\varepsilon)
&\lesssim\sup_{t\leq T}\int_{\mathbb R^3}
\left[|u|^\gamma-(|u|^2+\varepsilon^2)^{\gamma/2}\right]
\left(\int_{\mathbb R^3}|g(t,v_*+u)|\cdot|g(t,v_*)|\,dv_*\right)du\\
&\lesssim\left(\sup_{t\leq T}\norm{g(t)}_{L_v^2}^2\right)
\int_{|u|\leq1}\left[|u|^\gamma-(|u|^2+\varepsilon^2)^{\gamma/2}\right]du+
\frac{|\gamma|}{2}\varepsilon^2\sup_{s\leq T}\|g(s)\|_{L_v^1}^2\longrightarrow0,
\end{aligned}
\]
where we use the mean value theorem for $|u|\geq 1$.
In conclusion, discrete Gronwall yields, for $j\delta=t\leq T$,
\begin{equation}\label{eq:limitation for approximation}
l_j\leq te^{C_Tt}\left[C_\varepsilon\delta+\eta_1(\varepsilon)+\eta_2(\delta)\right].
\end{equation}

For each fixed $t>0$, choose $N_\varepsilon\in\mathbb N$ so that
$\delta_\varepsilon:=\dfrac{t}{N_\varepsilon}\leq\dfrac{\varepsilon}{1+C_\varepsilon}.$
Let $\varepsilon_n\downarrow0$. By
\eqref{eq:limitation for approximation},
$\|g_{N_{\varepsilon_n}}-g(t)\|_{L^1_v}\longrightarrow0.$
Since $g_{N_{\varepsilon_n}}\geq0$ and
$\int g_{N_{\varepsilon_n}}\,dv=M_0$, it follows that
$g(t)\geq0$ and $\int g(t)\,dv=M_0$. Consequently,
\[
\|\sqrt{g_{N_{\varepsilon_n}}}-\sqrt{g(t)}\|_{L^2_v}^2
\leq\|g_{N_{\varepsilon_n}}-g(t)\|_{L^1_v}
\longrightarrow0.
\]
Moreover, \eqref{eq:split-Fisher} and mass conservation give the boundedness of 
$\norm{\sqrt{g_{N_{\varepsilon_n}}}}_{H^1_v}^2$.
Choose a subsequence realizing the liminf of the Fisher
informations, and then a further subsequence converging weakly
in $H^1_v$. The strong $L^2_v$ convergence identifies its weak
limit as $\sqrt{g(t)}$. Hence, by weak lower semicontinuity,
\[
\mathcal I(g(t))
=4\|\nabla\sqrt{g(t)}\|_{L^2_v}^2
\leq\liminf_{n\to\infty}
   \mathcal I(g_{N_{\varepsilon_n}})
\leq e^t\mathcal I(g_0).
\]
Then we finish the proof of \eqref{eq:full-Fisher}.
\end{proof}

\subsection{Homogeneous data with finite Fisher information}
The proof of the homogeneous assertion uses only mass and Fisher
information. The approximation requires no finite energy, higher
velocity moments, or smallness of the initial data.

\begin{proof}[The proof of Theorem~\ref{thm:shear}
\textup{(\ref{part:shear-homogeneous})}]
\textit{Step 1: Local $Z$ solution and the blow up criterion}

Considering the ball
$\left\{u\in C([0,T];Z):\norm{u}_{C([0,T];Z)}\leq R\right\}$
and the mapping
$(\Phi u)(t)=S(t)g_0+\displaystyle\int_0^t S(t-s)\left[Q(u(s),u(s))\right]ds.$
By Lemma \ref{bi-esti-shear}, we deduce
\begin{align*}
    \norm{\Phi u}_{C([0,T];Z)} &\le \norm{g_0}_Z+CT R^2 \\
    \norm{\Phi u_1-\Phi u_2}_{C([0,T];Z)} &\le 2CT R\norm{u_1-u_2}_{C([0,T];Z)}.
\end{align*}
Omitting a standard contraction argument, we obtain a unique maximal solution
$g\in C([0,T_*);Z),$
where $T_*$ is the lifespan. For every $M>0$, every initial datum satisfying $g_0\geq0, \norm{g_0}_{Z}\leq M$ generates a local solution with lifespan at least $T(M)=\dfrac{1}{8C(1+M)}$.

Now we prove the blow up criterion as follows: if $ T_*<\infty$, then
    $$\lim_{T\uparrow T_*}\sup_{0\leq t\leq T}\norm{g(t)}_{Z}=\infty.$$
Just suppose that $T_*<\infty$, and there exists a sequence $t_j\uparrow T_*$ and a constant $M>0$ such that
$\|g(t_j)\|_{Z}\leq M.$
Restarting the local existence argument at each time $t_j$,
we extends the solution forward by at least $T(M)$. Choose $j$ sufficiently large that $T_*-t_j<T(M).$ The restarted solution then exists beyond $T_*$. A contradiction.

\medskip
\textit{Step 2: Fisher control and global existence.}

On every compact subinterval of $[0,T_*)$,
        Proposition~\ref{prop:Fisher} and Lemma~\ref{Endpoint-Loss term} give
\[
\norm{g(t)}_Z\lesssim M_0+\mathcal{I}(g(t))\leq M_0+e^t\mathcal I(g_0)<\infty
\]
By the blow up criterion in Step 1, $T_*=\infty$ and $g\in C([0,\infty);Z)$.
For every finite $T$, using Lemma~\ref{Endpoint-Loss term} again we have
\begin{equation}\label{eq:finite-time-bounds}
 \sup_{0\leq t\leq T}
 [\norm{g(t)}_Z+\norm{g(t)}_{L^3_v}+\norm{A_\gamma [g(t)]}_{L^\infty_v}]
 \lesssim M_0+e^T\mathcal{I}(g_0).
\end{equation}

\medskip
\textit{Step 3: Global uniqueness.}

Let $g_1,g_2\geq0$ be two mild solutions with the same initial datum $g_0$,
in $C([0,T];Z)$, and put $h=g_1-g_2$.
Then $\frac{d}{dt} [S(-t)h(t)]=S(-t)[Q(h,g_2)+Q(g_1,h)](t)$,
 which implies
\[
\norm{h(t)}_{L^1_v}\lesssim\int_0^t\bigl(\norm{A_\gamma[g_1(s)]}_{L^\infty_v}
+\norm{A_\gamma[g_2(s)]}_{L^\infty_v}\bigr)\norm{h(s)}_{L^1_v} ds.
\]
Combing with \eqref{eq:finite-time-bounds}, Gronwall yields the uniqueness.

\medskip
\textit{Step 4. Propagation of initially finite Sobolev norms.}
\

It's standard to obtain a unique solution $g\in C([0,T];Z^m)$ and, by zero-order uniqueness, coincides with the solution constructed in Step 2. Let $T_m^*$ denote the maximal existence time of $g$. Similarly, if $T_m^*<\infty$, then $\limsup_{t\uparrow T_m^*}\norm{g(t)}_{Z^m}=\infty.$

Now define
$\mathcal{L_{H}}=\partial_t-v_2\partial_{v_1},\ \tilde{Y}_1=\partial_{v_1},\ \tilde{Y}_2=\partial_{v_2}-t\partial_{v_1},\ \tilde{Y}_3=\partial_{v_3}.$
A direct calculation shows:
$[\mathcal{L_{H}},\tilde{Y}_j]=0, [\tilde{Y}_i,\tilde{Y}_j]=0.$
For a multi-index $\beta=(\beta_1,\beta_2,\beta_3)$, considering
$$
\tilde{Y}^\beta(t)g(t)=S(t)\partial_v^\beta g_0
+\sum_{\eta\leq\beta}\binom{\beta}{\eta}\int_0^tS(t-s) Q\bigl(\tilde{Y}^\eta(s)g(s), \tilde{Y}^{\beta-\eta}(s)g(s)\bigr)\,ds.
$$
Define $N_j(t)=\sum_{|\beta|=j}\norm{\tilde{Y}^\beta(t)g(t)}_{Z},\
\mathfrak N_j(T)=\sup_{0\leq t\leq T}N_j(t),$ then
$$
N_m(t)\lesssim N_m(0)+ \int_0^t \left[ N_0(s)N_m(s)+\sum_{j=1}^{m-1}N_j(s)N_{m-j}(s) \right]ds.
$$
Now begin induction, for $m=1$, Gronwall's inequality gives
$\mathfrak N_1(T)\lesssim N_1(0)e^{C_m\mathfrak N_0(T)T}<\infty.$
If $\mathfrak N_j(T)<\infty$ for every $j<m$, then
$$
\mathfrak N_m(T)\lesssim\left[ N_m(0)+T\sum_{j=1}^{m-1}\mathfrak N_j(T)\mathfrak N_{m-j}(T)\right] e^{C_m\mathfrak N_0(T)T}<\infty.
$$
Since
$\widetilde Y_j(t)=S(t)\partial_{v_j}S(-t)$
and $S(t)$ is an isometry on $Z$, we have
\[
\|S(-t)g(t)\|_{Z^m}
\sim \sum_{j=0}^m N_j(t), \qquad \|S(t)h\|_{Z^m}\leq C_m(1+|t|)^m\|h\|_{Z^m}
\]
Consequently, for $0\leq t\leq T<T_m^*$,
\[
\|g(t)\|_{Z^m}
\leq C_m(1+T)^m\sum_{j=0}^m N_j(t)<\infty,
\]
which,by the blow-up criterion, implies
$T_m^*=\infty,$ and $ g\in C([0,\infty);Z^m).$
\end{proof}
The preceding theorem also applies to sufficiently regular initial data, since their Fisher information can be controlled by Sobolev norms.
\begin{corollary}
Let $m\geq3$ be an integer and $-2\leq\gamma\leq0$.
If $g_0\geq0$ belongs to $Z^m$, then
Theorem~\ref{thm:shear}\textup{(\ref{part:shear-homogeneous})}
holds without separately assuming $\mathcal I(g_0)<\infty$.
In particular, the unique nonnegative global solution satisfies
$g\in C([0,\infty);Z^m)$ and
\begin{equation}\label{eq:high-regularity-fisher}
 \mathcal I(g(t)) \leq e^t\mathcal I(g_0)
 \lesssim e^t M_0^{1/3}\left(\sum_{i=1}^3\|\partial_{v_i}^3g_0\|_{L^1_v}\right)^{2/3}.
\end{equation}
\end{corollary}

\begin{proof}
Without loss of generality assume $M_0>0$, and let $g_{0,\eta}=\rho_\eta*g_0$, where $\rho_\eta\geq 0$ is a standard mollifier. Applying \cite[Theorem~1.2]{Bobkov2022} to the probability density $g_{0,\eta}/M_0$ gives
\[
 \mathcal I(g_{0,\eta})\lesssim
 M_0^{1/3}\left(\sum_{i=1}^3\|\partial_{v_i}^3g_{0,\eta}\|_{L^1_v}\right)^{2/3}
 \lesssim M_0^{1/3}\left(\sum_{i=1}^3\|\partial_{v_i}^3g_0\|_{L^1_v}\right)^{2/3}.
\]
Since $g_{0,\eta}\to g_0$ in $L^1_v$, lower semicontinuity of
Fisher information yields \eqref{eq:high-regularity-fisher} and thus
the first assertion.
\end{proof}

Finally, the Fisher bound yields decay estimates after a mass-preserving dilation, as described in the following remark.
\begin{remark}
The Fisher bound has a simple consequence under mass-preserving dilation. For $a>1/2$, set $F(t,\xi)=e^{-3at}f(t,e^{-at}\xi)$. Direct substitution gives
\begin{equation*}
 \partial_tF-\xi_2\partial_{\xi_1}F +a\nabla_\xi\cdot(\xi F)=e^{-\gamma at}Q(F,F).
\end{equation*}
The change of variables preserves mass and gives
\begin{equation*}
 \mathcal I(F(t))=e^{-2at}\mathcal I(f(t)) \leq e^{-(2a-1)t}I_0.
\end{equation*}
Thus, for $1\leq p\leq3$,
\[
 \norm{F(t)}_{L^p_\xi} \lesssim e^{-(2a-1)(1-\frac1p) t}.
\]

\end{remark}

\subsection{Small spatially inhomogeneous data}
We now close the global estimate for \eqref{eq:shear-cauchy} using the
vector fields and bilinear estimates established in Section~\ref{dis-esti}.
For the simple shear flow, the matrices in \eqref{eq:commuting-matrices}
take the explicit form
\[
 \Phi(t)=\begin{pmatrix}1&-t&0\\0&1&0\\0&0&1\end{pmatrix},
 \qquad
 \Psi(t)=\begin{pmatrix}1+t&-t^2/2&0\\0&1+t&0\\0&0&1+t\end{pmatrix}.
\]
In particular, $\det \Phi(t)=1$ and $\det \Psi(t)=(1+t)^3$.
Applying Lemma~\ref{thm:dis-esti}, we therefore obtain, for
$1\leq p<\infty$,
\begin{equation}\label{eq:shear-dispersion}
 \|h(t)\|_{L_x^{\infty} L_v^p}
 \lesssim_p (1+t)^{-3/p}
 \sum_{|I|\leq3}\|\mathcal D^I(t)h(t)\|_{L_{x,v}^p}.
\end{equation}
By a standard density argument, this estimate extends from smooth
functions to the corresponding function spaces. Consequently,
combining \eqref{eq:shear-dispersion} with the collision estimates
established in Section~\ref{dis-esti} yields the bilinear estimate
\eqref{Estimate-Q-X}. 

\begin{proof}[Proof of Theorem~\ref{thm:shear}\textup{(\ref{part:shear-inhomogeneous})}]
Let $U(t)$ be the solution group for
$\partial_t+v\cdot\nabla_x-v_2\partial_{v_1}$; explicitly,
\[
 (U(t)h)(x,v)
 =h\bigl(x-tv-\tfrac12t^2v_2e_1,\,v+tv_2e_1\bigr).
\]
It is an isometry on every $L^p_{x,v}$, and
$\mathcal D^I(t)U(t-s)=U(t-s)\mathcal D^I(s)$.
Consequently, the mild equation and the bilinear estimate
\eqref{Estimate-Q-X} yield
\begin{equation}\label{priori-shear-1}
\begin{aligned}
     \sum_{|I|\leq N_p}\norm{\mathcal{D}^I f}_{L_{x,v}^p}\lesssim& \sum_{|I|\leq N_p}\norm{\mathcal{D}^{I}(0)f_0}_{L_{x,v}^p}+\sum_{|I|\leq N_p}\int_{0}^{T}\norm{\mathcal{D}^{I}Q(f,f)}_{L_{x,v}^p}ds\\
     \lesssim&E_0(f_0) + \int_{0}^{T}(1+s)^{-3/p_*}[E_{s}(f)]^{2}ds\\
     \lesssim&E_0(f_0)+[E_{T}(f)]^{2},
\end{aligned}
\end{equation}
where $N_p = \left\{
\begin{array}{cc}
    6, &  p=1    ;\\
   8,  &  p=p_+  ;\\
   11,  & p=p_{-},\\
\end{array}
\right.$
then \eqref{priori-shear-1} implies
\begin{equation}\label{priori-shear}
    E_T(f)\lesssim E_0(f_0)+[E_{T}(f)]^{2}.
\end{equation}
Therefore a priori estimate \eqref{priori-shear} and standard contraction argument  with sufficiently small initial value, bootstrap argument allow us  to extend this local solution globally in time, with
$$E_{\infty}(f) \lesssim E_0(f_0).$$
Moreover, the solution map $f_0\mapsto f$ is Lipschitz continuous, and the uniqueness follows
the difference estimate
$$E_T(f_1-f_2)\lesssim E_0(f_{1,0} - f_{2,0}) +[E_T(f_1)+E_{T}(f_2)]E_{T}(f_1 - f_2).$$
If $f_0\geq0$, the gain--loss formula along the characteristics,
first for the local iterates and then by convergence, preserves
nonnegativity, we omit the standard details. This ends the proof of theorem~\ref{thm:shear}\textup{(\ref{part:shear-inhomogeneous})},  including
the Maxwell endpoint.
\end{proof}

\appendix
\section{Auxiliary inequalities}\label{sec:appendix}
This appendix collects the elementary Sobolev and interpolation
inequalities and the velocity-space collision bounds used in
Section~\ref{dis-esti}. We state the exponent restrictions explicitly,
since the gain and loss terms have different endpoint behavior.

\begin{lemma}
\label{lem:banach-valued-sobolev}
Let $V$ be a Banach space and $1\leq p<\infty$. If
$u\in W^{3,p}(\mathbb R^3;V)$, then $u$ has a bounded representative and
\begin{equation*}
 \|u\|_{L_x^\infty V}
 \leq C_p\sum_{|\alpha|\leq3}
             \|\partial_x^\alpha u\|_{L_x^pV}.
\end{equation*}
In particular, the assertion applies to $V=L_\xi^q$ for
$1\leq q<\infty$.
\end{lemma}
The proof of the above proposition can be found in \cite[Theorem 6.1]{WM}.



\begin{lemma}
\label{lem:collision-velocity}
      Let $1<p,q,r<\infty$ and $\gamma\in (-d,0]$ with $$ \dfrac{1}{p}+\dfrac{1}{q} = 1 + \dfrac{\gamma}{d} + \dfrac{1}{r}.$$
     Assume that the collision kernel satisfies the angular cutoff condition \eqref{eq:angular-cutoff}. Then the following estimates hold:
      \begin{align*}
        \norm{Q^{+}(f,g)}_{L_v^r(\mathbb{R}^d)}\lesssim \norm{f}_{L_v^p(\mathbb{R}^d)}\norm{g}_{L_v^q(\mathbb{R}^d)}, &  \\
        \norm{Q^{-}(f,g)}_{L_v^r(\mathbb{R}^d)}\lesssim \norm{f}_{L_v^p(\mathbb{R}^d)}\norm{g}_{L_v^q(\mathbb{R}^d)}, &\quad\quad p>r.
      \end{align*}
    \end{lemma}
    The details of the lemma above can be found in
\cite[Lemma A.4]{ChenShenZhangSoft}. 

For later use, set
\[
 A_\gamma[h](v)=\int_{\mathbb R^3}|v-v_*|^\gamma h(v_*)\,dv_*.
\]
\begin{lemma}
\label{lem:collision-frequency}
Let $-3<\gamma<0$, $p_c=3/(3+\gamma)$, and
$1\leq r<p_c<q\leq\infty$. Then
\begin{equation*}
 \big\|A_\gamma[h]\big\|_{L_v^\infty}
 \lesssim\|h\|_{L_v^r}+\|h\|_{L_v^q}.
\end{equation*}
\end{lemma}
\begin{proof}
We split the kernel into the regions $|v|\leq 1$ and $|v|>1$. By H\"{o}lder's inequality,
\begin{align}
\norm{A_\gamma[h]}_{L_v^{\infty}}
=\norm{|\cdot|^\gamma*h}_{L_v^{\infty}}
\leq
\norm{|v|^\gamma\mathbf 1_{{|v|\leq1}}}_{L_v^{q'}}
\norm{h}_{L_v^{q}}
+
\norm{|v|^\gamma\mathbf 1_{{|v|>1}}}_{L_v^{r'}}
\norm{h}_{L_v^{r}}.\notag
\end{align}
The first kernel norm is finite since
$
\gamma q'>-3.
$
If $r>1$, the second kernel norm is finite provided
$
\gamma r'<-3.
$
If $r=1$, then $r'=\infty$, and
$
\norm{|v|^\gamma\mathbf 1_{{|v|>1}}}_{L_v^\infty}\leq 1,
$
since $\gamma<0$. The desired estimate follows.
\end{proof}

\begin{lemma}
\label{Endpoint-Loss term}
    Let $-2\leq\gamma\leq 0$.
    \begin{enumerate}
        \item [(1)] If $h\in L^1,\ h\geq 0$ and $\mathcal{I}(h)
    <\infty$, then
    \begin{equation}\label{eq:Fisher-Lorentz}
    \norm{h}_{L_{v}^{3,1}}+\norm{h}_{L_v^3}\lesssim \mathcal{I}(h).
    \end{equation}
        \item [(2)] If $h\in Z$,
    then
    \begin{equation}\label{eq:Z norm-Loss term}
    \norm{A_{\gamma}[|h|]}_{L_v^\infty}\lesssim \norm{h}_Z.
    \end{equation}
    \end{enumerate}
\end{lemma}

\begin{proof}
We first prove \eqref{eq:Fisher-Lorentz}. In view of the Sobolev embedding
$H^1\hookrightarrow L^6$, it suffices to show that
$\norm{h}_{L_v^{3,1}} \lesssim \mathcal{I}(h).$

Set $u=\sqrt{h}$. By Hardy's inequality, we have
\[
\norm{h}_{L_v^{3,1}}
\sim \int_{\mathbb{R}^3}\frac{[u^*(|v|^3)]^2}{|v|^2} dv
\lesssim \int_{\mathbb{R}^3}\moc{\nabla u^*(|v|^3)}^2 dv.
\]
It remains to control the last term by the Dirichlet energy of $u$.
To this end, let
$E_{\lambda}=\{v\in\mathbb{R}^3: u(v)>\lambda\}$.
By the Cauchy--Schwarz inequality, the coarea formula, and the
isoperimetric inequality \cite[Theorem~2.17]{Burchard}, we obtain
\[
\int_{\{u=\lambda\}} |\nabla u| dS
\geq
\frac{|\partial E_\lambda|^2}{\int_{\{u=\lambda\}} 1/|\nabla u| dS}
\gtrsim
\frac{|E_\lambda|^{\frac 43}}{\int_{\{u=\lambda\}} 1/|\nabla u| dS}.
\]
On the other hand, applying the same argument to the radial rearrangement
$u^*(|v|^3)$ gives
\[
\int_{\{ u^*(|v|^3)=\lambda \}} | \nabla u^*(|v|^3) | dS
\sim
\frac{|E_\lambda|^{\frac 43}}
{\int_{\{ u^*(|v|^3)=\lambda \}} 1/| \nabla u^*(|v|^3) | dS}
\sim
\frac{|E_\lambda|^{\frac 43}}{\int_{\{u=\lambda\}} 1/|\nabla u| dS},
\]
where we use
$|\partial \{ u^*(|v|^3)>\lambda \}|=4\pi |E_\lambda|^{\frac 23}$.
Combining the above estimates and applying the coarea formula once more,
we conclude that
\[
\int_{\mathbb{R}^3}[\nabla u^*(|v|^3)]^2 dv
=\int_0^\infty d\lambda\int\limits_{\{ u^*(|v|^3)=\lambda \}}
| \nabla u^*(|v|^3) | dS\\
\lesssim \int_0^\infty d\lambda \int\limits_{\{u=\lambda\}} |\nabla u| dS
= \int_{\mathbb{R}^3} |\nabla u|^2 dv,
\]
which proves \eqref{eq:Fisher-Lorentz} by a standard approximation
argument.

We next prove \eqref{eq:Z norm-Loss term}. Without loss of generality,
assume that $\gamma\in [-2,0)$. Since
$|v|^{-2}\in L_v^{\frac32,\infty}$, the Lorentz H\"older inequality
$\int FG\leq \norm{F}_{L^{3,1}}\norm{G}_{ L^{\frac32,\infty}}$
implies
$\norm{A_{-2}[h]}_{\displaystyle L_v^\infty}
\lesssim \norm{h}_{\displaystyle L_v^{3,1}}.$
Splitting the integral into the regions $|v-v_*|\leq1$ and
$|v-v_*|\geq1$, we therefore obtain
\[
\norm{A_\gamma[h]}_{L_v^\infty}\leq
\norm{\int_{|v-v_*|\leq 1} |v-v_*|^{-2} h dv_*}_{L_v^\infty}+
\norm{\int_{|v-v_*|\geq 1} 1\cdot h dv_*}_{L_v^\infty}
\lesssim \norm{h}_Z.
\]
This proves \eqref{eq:Z norm-Loss term} and completes the proof.
\end{proof}

\noindent {\bf Acknowledgements:}
	SQL was supported by grants from the National Natural Science Foundation of China under contracts 12325107 and 12471217.

\medskip
\noindent\textbf{Data Availability Statement:}
Data sharing is not applicable to this article as no datasets were generated or analysed during the current study.

\noindent\textbf{Conflict of Interest:}
The authors declare that they have no conflict of interest.

\noindent\textbf{AI usage disclosure:}  The authors used AI-assisted tools, principally ChatGPT 5.6, for writing assistance and to help check the correctness of the mathematical arguments. The authors have carefully reviewed all AI-assisted revisions and take full responsibility for the content, arguments, and conclusions presented in the manuscript.

\end{document}